\documentclass{amsart}

  \usepackage{graphicx}
  \usepackage{amsfonts}
  \usepackage{amsmath}
  \usepackage{amsthm,amstext,amsfonts,bm,amssymb}
  \usepackage[normalem]{ulem}
  \usepackage[below]{placeins}
  \usepackage{amssymb,enumerate,hyperref}
  \usepackage{xfrac}
  \usepackage{tikz-cd}
  \usepackage{refcount}
  \usepackage{hyperref} 
  \usepackage[nameinlink]{cleveref}
  \usepackage[margin=1cm,singlelinecheck=off]{caption}
  \usepackage{enumitem,mathtools}

  \providecommand{\U}[1]{\protect\rule{.1in}{.1in}}
  \providecommand{\comment}[1]{}
  \newcommand{\actson}{\curvearrowright}
  		\newcommand{\CB}{\mathcal B}
  \newcommand{\CE}{\mathcal E}		

  		\newcommand{\CN}{\mathcal N}

  \newcommand{\CS}{\mathcal S}

  \renewcommand{\int}{\text{int}}
  \providecommand{\co}{\mskip0.5mu\colon\thinspace} 

  \newcommand{\BN}{\mathbb N}
  \newtheorem{theorem}{Theorem}[section]

  \newtheorem*{namedtheorem}{\theoremname}
  \newcommand{\theoremname}{testing}
  \newenvironment{named}[1]{\renewcommand{\theoremname}{#1}\begin{namedtheorem}}{\end{namedtheorem}}

  \makeatletter
  \newtheorem*{rep@theorem}{\rep@title}
  \newcommand{\newreptheorem}[2]{%
  \newenvironment{rep#1}[1]{%
   \def\rep@title{#2 \ref{##1}}%
   \begin{rep@theorem}}%
   {\end{rep@theorem}}}
  \makeatother

  \newreptheorem{theorem}{Theorem}
  \newreptheorem{lemma}{Lemma}
  \newreptheorem{proposition}{Proposition}
  \newtheorem{claim}[theorem]{Claim}

  \newtheorem{corollary}[theorem]{Corollary}

  \newtheorem{lemma}[theorem]{Lemma}
  
  \newtheorem*{lemma*}{Lemma}

  \newtheorem{proposition}[theorem]{Proposition}

  \theoremstyle{definition}
    \newtheorem{remark}[theorem]{Remark}
  \newtheorem{definition}[theorem]{Definition}
  \newtheorem{question}[theorem]{Question}
  
  \newtheorem*{definition*}{Definition}

  \newcommand{\BZ}{\mathbb{Z}}
  
  \newcommand{\BR}{\mathbb{R}}

  \newcommand{\dsm}{\mathbin{{\smallsetminus}\mspace{-5mu}{\smallsetminus}}}

\newcommand{\BCT}{\mathcal B}

  \usepackage{todonotes}

	\definecolor{cadmiumorange}{rgb}{0.93, 0.53, 0.18}

  \title{Covers of surfaces}

 \author[Biringer]{Ian Biringer}
 \address{\parbox{\linewidth}{Ian Biringer, Boston College, 140 Commonwealth Ave, Chestnut Hill, MA 02467\vspace{1mm}}}
\email{ian.biringer@bc.edu}
  
  \author[Chandran]{Yassin Chandran} 
 \address{\parbox[t]{\linewidth}{Yassin Chandran, CUNY Graduate Center, 365 Fifth Avenue, New York, NY 10016\vspace{1mm}}}
\email{ychandran@gradcenter.cuny.edu}
 
  \author[Cremaschi]{Tommaso Cremaschi} 
\address{\parbox{.95\linewidth}{Tommaso Cremaschi, School of Mathematics, Trinity College Dublin, 17 Westland Row, Dublin 2, Ireland\vspace{1mm}}}
\email{cremasct@tcd.ie}
  
  \author[Tao]{Jing Tao}
  \address{\parbox{.95\linewidth}{Jing Tao, Department of Mathematics, University of Oklahoma,
  601 Elm Ave, Room 423, Norman, OK 73069
  \vspace{1mm}}}
\email{jing@ou.edu}
  
  \author[Vlamis]{Nicholas G.~Vlamis}
  \address{\parbox{.95\linewidth}{Nicholas G.~Vlamis, Department of Mathematics, CUNY Graduate Center, 365 Fifth Avenue, New York, NY 10016, and 
  Department of Mathematics, CUNY Queens College, 65-30 Kissena Blvd., Flushing, NY 11367\vspace{1mm}}}
  \email{nvlamis@gc.cuny.edu}
  
  \author[Wang]{Mujie Wang}
   \address{\parbox{\linewidth}{Mujie Wang, Boston College, 140 Commonwealth Ave, Chestnut Hill, MA 02467\vspace{1mm}}}
  \email{mujiew@bc.edu}
   
  \author[Whitfield]{Brandis Whitfield}
  \address{\parbox{.95\linewidth}{Brandis Whitfield, University of Wisconsin, 480 Lincoln Dr., Madison, WI 53706\vspace{1mm}}}
  \email{brandis.whitfield@wisc.edu}

\date{
    \today
}

\begin{document}

  \thispagestyle{empty}

  \begin{abstract} 
We study the homeomorphism types of certain covers of (always orientable) surfaces, usually of infinite type. We show that every surface with non-abelian fundamental group is covered by every noncompact surface, we identify the universal abelian covers and the $\BZ/n\BZ$-homology covers of surfaces, and we show that non-locally finite characteristic covers of surfaces have four possible homeomorphism types. 
  \end{abstract}

\maketitle

\tableofcontents
  
\section{Introduction}

We are interested in the question of when a given surface covers another given surface, either in general or with constraints on the type of covering. As motivation, let us mention the following elementary result.

  \begin{proposition}[Everything covers everything]\label{Prop:arbitrary covers}
    Suppose that $S$ is an orientable, borderless surface with non-abelian fundamental group. Then $S$ is covered by any noncompact borderless orientable surface.
 \end{proposition}

A \emph{borderless} surface is a topological $2$-manifold with empty boundary. We stick with orientable surfaces for simplicity. The only orientable surfaces $S$ not covered by the theorem are the disk, the annulus, and the torus, each of which is only covered by itself and the former surface. 

Proposition \ref{Prop:arbitrary covers} is not surprising, and it is not difficult. Related results already appear in the literature: Goldman \cite{Goldman1971} proves the result for closed $S$, which may also be non-orientable, and similar results are implicit in some papers on foliations, see e.g.\ \cite{CC87}, \cite{AB22}, \cite{ABPW22}, and \cite{GM24}. With some extra work, the methods used by Goldman and others could likely yield our statement. However, our proof is short, constructive and the result motivates well the rest of our work.

\medskip

In another direction, every surface $S$ has distinguished covers coming from especially important subgroups of $\pi_1 (S)$. As a simple first example, the \emph{universal abelian cover} of $S$ is the cover $S^{ab} \to S$ corresponding to the commutator subgroup of $\pi_1(S)$. In contrast to the flexibility displayed in Proposition \ref{Prop:arbitrary covers}, universal abelian covers only have five possible homeomorphism types!

  \begin{theorem}[Universal abelian covers]\label{Prop:UAC}
    Let $S$ be an orientable, borderless surface, and let $S^{ab}$ be its universal abelian cover. Then
    \begin{enumerate}
      \item if $S$ is $\BR^2$, the annulus, or the torus then $S^{ab}\cong \BR^2$,
      \item if $S$ is the sphere, so is $S^{ab}$,
      \item if $S$ is the once-punctured torus, then $S^{ab}$ is a flute surface,
      \item if $S$ is a finite-type surface with one puncture and genus at least $2$, then $S^{ab}$
        is the spotted Loch Ness monster surface,
      \item otherwise, $S^{ab}$ is the Loch Ness monster surface. 
    \end{enumerate}
  \end{theorem}

The proof uses some basic geometric group theory to identify the universal abelian covers of finite type surfaces, and then for the infinite type case, we exhaust by finite type subsurfaces.

The commutator subgroup of $\pi_1(S)$ is the kernel of the map $\pi_1(S) \to H_1(S,\mathbb Z),$ so the universal abelian cover could be called the {``$\BZ$-homology cover''} of $S$. We also consider the \emph{$\BZ/n\BZ$-homology cover} $\tilde S \to S$ corresponding to the kernel of $$\pi_1(S) \to H_1(S,\BZ/n\BZ).$$ When $S$ has finite-type, the $\BZ/n\BZ$-homology cover $\tilde S \to S$ is a finite cover, and one can figure out the topological type of $\tilde S$ via an  Euler characteristic argument and by analyzing how boundary components lift. On the other hand, when $\tilde S$ has infinite-type, its $\BZ/n\BZ$-homology covers all have infinite degree, and it turns out their topology is constrained, just like for $\BZ$-homology covers.

\begin{theorem}[$\BZ/n\BZ$-homology covers]\label{thm: finite homology covers}
    Suppose that $S$ is an infinite-type orientable borderless surface, $n\geq 2$, and $\pi: \tilde S \to S$ is the $\BZ/n\BZ$-homology cover.
    If $S$ has no isolated planar ends, then $\tilde S $ is homeomorphic to the Loch Ness monster surface. Otherwise, $\tilde S$ is homeomorphic to the spotted Loch Ness monster surface.
\end{theorem}

Suppose that a cover $\tilde S \to S$ corresponds to the subgroup $N < \pi_1 (S,p)$. We say that the cover is \emph{characteristic} if $N$ is characteristic, i.e.,\ if $N$ is invariant under every automorphism of $\pi_1(S,p)$. The $\BZ$ and $\BZ/n\BZ$-homology covers above are all characteristic. Motivated by Proposition \ref{Prop:UAC} and Theorem \ref{thm: finite homology covers}, we ask:

\begin{question}\label{charquest}
    Are all infinite-degree characteristic covers of orientable borderless surfaces homeomorphic to either the disk, the flute surface, the Loch Ness monster surface, or the spotted Loch Ness monster surface?
\end{question}

Note that it is important to say `infinite degree' above. Indeed, if $S$ has finite type it has many finite-degree characteristic covers $\tilde S$, all of which also have finite type. (For these $\tilde S$, in practice one can usually figure out the topological type of $\tilde S$ by analyzing how boundary components lift and using an Euler characteristic argument.) When $S$ has infinite type, the `infinite degree' hypothesis just prohibits the cover from being trivial (which is obviously necessary), since $F_\infty$, the free group on a countably infinite alphabet, does not have any proper finite-index characteristic subgroups.

Here is some evidence that the answer to Question \ref{charquest} could be yes.

\begin{theorem}[Characteristic covers of finite-type surfaces]\label{thm: finite type}
    Suppose that $S$ is an orientable borderless surface of finite type and $\tilde S \to S$ is an infinite-degree, geometrically characteristic cover. Then $\tilde S $ is homeomorphic to either a disc, the flute surface, the Loch Ness monster surface, or the spotted Loch Ness monster surface.
\end{theorem}

An automorphism of $\pi_1(S,p)$ is called \emph{geometric} if it is induced by a homeomorphism $S \to S$ that preserves the basepoint $p$, a subgroup $N < \pi_1(S,p)$ is \emph{geometrically characteristic} if $N$ is invariant under all geometric automorphisms of $\pi_1 (S,p)$. A cover $\tilde S \to S$ is \emph{geometrically characteristic} if it corresponds to a geometrically characterstic subgroup of $\pi_1(S,p)$. Of course, any characteristic cover is geometrically characteristic. Note that geometrically characteristic subgroups are normal, since inner automorphisms of $\pi_1(S,p)$ can be realized via point pushing homeomorphisms, and hence geometrically characteristic subgroups are regular. Geometrically characteristic covers were previously considered in \cite{aramayona2023big}, where they were defined equivalently as regular covers to which every homeomorphism lifts\footnote{There are non-regular covers of closed surfaces where every homeomorphism lifts, see Propositions 6 and 7 in \cite{aramayona2023big}. These are not geometrically characteristic as defined above.}. 

It follows quickly from Theorem \ref{thm: finite type} that any non-universal characteristic cover $\tilde S$ of a \emph{closed} surface is homeomorphic to the Loch Ness monster surface. This was proved earlier by Atarihuana--Garc\'ia--Hidalgo--Quispe--Maluendas \cite[Prop 6.1]{atarihuana2022dessins}, using the fact that $\tilde S$ is also a characteristic cover of a triangle orbifold. However, their argument does not extend to the finite-type setting.

Infinite-type surfaces can be exhausted by finite-type surfaces, and characteristic subgroups of the infinite-type surface group intersect the finite-type surface groups in characteristic subgroups (see \S \ref{sec:pfchar} for details) which suggests that the theorem above can be promoted to apply to infinite-type $S$. However, there is a problem: there are infinite-degree characteristic covers of infinite-type surfaces $\tilde S$ that restrict to finite-degree covers on all finite-type subsurfaces, and Theorem~\ref{thm: finite type} does not apply to finite covers. To deal with this, we make the following definition.

\begin{definition}\label{locallyfinite}
    If $G$ is a group, a subgroup $N < G$ is \emph{locally finite index} if $N$ intersects every finitely generated subgroup $F \subset G$ in a finite-index subgroup of $F$. A cover $\pi: \tilde S \to S$ is \emph{locally finite} if it corresponds to a locally finite index subgroup of $\pi_1(S) \cong F_\infty$. Equivalently, $\pi$ is locally finite if for every compact subsurface $K\subset S$, each component of $\pi^{-1}(K)$ is compact.
\end{definition}

When $S$ has finite-type, a locally finite cover $\tilde S \to S$ is just a finite cover. When $S$ has infinite-type, the $\BZ/n\BZ$-homology covers discussed above are locally finite, but not finite. Locally finite index subgroups have been studied previously in  \cite{banakh2010coarse,hampton1972semisimplicity}, for instance, but we are not aware of any previous discussion of locally finite covering maps. We prove:

\begin{theorem}[Characteristic covers that are not locally finite]\label{thm: not locally finite}
        Suppose that $S$ is an orientable, infinite-type borderless surface and $\tilde S \to S$ is a characteristic cover that is not locally finite. Then $\tilde S $ is homeomorphic to either a disc, the flute surface, the Loch Ness monster surface, or the spotted Loch Ness monster surface.
\end{theorem}

So, to answer Question \ref{charquest}, it only remains to understand locally finite characteristic covers of infinite-type surfaces. One can check that the $\BZ/n\BZ$-homology covers discussed above are the only such covers that are abelian, and we have answered  Question~\ref{charquest} positively for those covers. An example of a non-abelian, locally finite, characteristic cover of a surface $S$ is the cover $\tilde S$ corresponding to the subgroup $N<\pi_1 (S)$ generated by all third powers of elements of $\pi_1(S)$. This cover is obviously characteristic, and it is locally finite because Burnside \cite{burnside1902unsettled} proved that the group $B(d,3) = F_d / \langle \langle w^3 \ | \ w\in F_d \rangle \rangle $ is finite for all $d$. It is not clear to us what $\tilde S$ is in this case.

Our interest in Theorem \ref{thm: not locally finite} is partly inspired by work on big mapping class groups of surfaces. In particular, Aramayona--Leininger--McLeay  \cite[Proposition~3.3]{aramayona2023big} show that the $\BZ/2\BZ$ homology cover of the Loch Ness monster surface is the Loch Ness monster surface. They use this to show that the mapping class group $G$ of the once-punctured Loch Ness monster surface admits a continuous injective homomorphism $G\to G$ that is not surjective. We note that there is a gap in the published version of their proof of the covering statement, but they were able to fix it quickly when we notified them about it. Alternatively, their Proposition 3.3 follows from 
Theorem \ref{thm: finite homology covers}, or even from Theorem \ref{thm: finite type} or  \cite[Proposition 6.1]{atarihuana2022dessins}, since the Loch Ness monster surface is itself a characteristic cover of a genus two surface.

Surprisingly, to our knowledge there are not many previous results in which the homeomorphism types of infinite covers of surfaces are studied.
In particular, we were surprised that a result like Proposition~\ref{Prop:arbitrary covers} does not currently exist in the literature. However, there are a few results out there that are related to ours. 
First, a classification of the possible homeomorphism types of infinite-degree regular covers of a closed surface is readily deduced from the fact that the cover and the deck group are quasi-isometric (for instance, see \cite[Proposition~5.2]{BasmajianExotic}).
Then there are older papers of Maskit \cite{maskit1965theorem} and Papakyriakopoulos \cite{papakyriakopoulos1975planar} on identifying when certain covers of surfaces are planar; their motivation was some interesting applications to $3$-manifolds, including a new proof of the Loop Theorem. There are the papers of Aramayona--Leininger--McLeay  \cite{aramayona2023big} and Atarihuana et al \cite{atarihuana2022dessins} referenced above, and there is also a forthcoming paper of Ghaswala--McLeay in which (among other things) the authors construct finite-degree covers from certain maps between potential end spaces of surfaces.

\subsection{Ideas of the proofs}

The proofs of our theorems above have distinct flavors. We will sketch each of them here for the convenience of the reader.

\medskip

Proposition \ref{Prop:arbitrary covers}, in which we show that every noncompact surface covers every surface with non-abelian fundamental group, is very hands-on. Every surface $S$ with non-abelian fundamental group has an essential pair of pants and so is covered by an open pair of pants. The universal abelian cover of an open pair of pants is the Loch Ness monster surface, which is covered by the blooming Cantor tree surface, and every noncompact surface $\tilde S$ can be $\pi_1$-injectively embedded in the blooming Cantor tree surface, and hence covers the blooming Cantor tree surface. (Here the embedding is as a component of the complement of a set of closed curves and bi-infinite arcs.) Composing, we get a covering map $\tilde S \to S$.

\medskip

Proposition \ref{Prop:UAC}, in which we identify the universal abelian covers $S^{ab}$ of surfaces $S$, is proved as follows. When $S$ has finite-type, one can identify $S^{ab}$ using the fact that the deck group is free abelian, along with the knowledge of whether the peripheral curves of $S$ lift or not. We deal with infinite-type $S$ by exhausting them with finite-type subsurfaces and saying that $S^{ab}$ is exhausted by the universal abelian covers of the subsurfaces.

\medskip

Theorem \ref{thm: finite homology covers}, in which we show that the $\BZ/n\BZ$-homology cover $\tilde S \to S$ of any infinite-type surface is either the Loch Ness monster surface or its spotted version, uses the fact that the cover is abelian, together with the fact that it unwraps curves that appear basically everywhere in $S$. Assuming for simplicity that $S $ has no cusps, we show that $\tilde S$ has one end as follows: take two points $\tilde p,\tilde q\in \tilde S$ that are far away from some chosen basepoint, and a path $\tilde a$ between them. This $\tilde a $ may pass close to the basepoint of $\tilde S$, but we can construct a new path $\tilde \gamma$ between them that stays far from the basepoint, by using abelianness of the cover to form a path-quadrilateral in $\tilde S$ where $\tilde a,\tilde \gamma$ are a pair of opposite sides. See also the beginning of \S \ref{sec: pflem}  for a more detailed sketch of the most important part of the argument.

\medskip

Theorem \ref{thm: finite type}, in which we identify all infinite-degree characteristic covers $\tilde S$ of finite-type surfaces $S$, has the most complicated proof. The argument is inspired by Proposition 3.3 in \cite{aramayona2023big}, which we discussed above.

    Say for simplicity that $S$ is closed. We want to prove that $\tilde S$ has one end. Here, $\tilde S$ is quasi-isometric to the deck group $Q$. If $Q$ has more than one end, it admits an action on a tree $T$ with finite edge stabilizers and no global fixed point, by Stallings's Theorem. We then construct a $Q$-equivariant map $\tilde S \to T$ that is transverse to the union of edge midpoints, and let $\tilde C \subset \tilde S$ be the preimage of the set of edge midpoints, which is then a properly embedded $1$-submanifold invariant under $Q$. Since edge stabilizers are finite, each component of $\tilde C$ is a simple closed curve. Using similar methods and an inductive argument, we can enlarge $\tilde C$ until every component of $\tilde S \smallsetminus \tilde C$ is one-ended. 
    
    If $C \subset S$ is the projection of $\tilde C$, and $\gamma$ is an essential curve in $S \smallsetminus C$, then $\gamma$ is `at most one ended', meaning that all lifts of $\gamma$ to $\tilde S $ are either compact or are bi-infinite arcs both of whose ends go out the same end of $\tilde S$. Any at most one-ended $\gamma$  projects to an element of $Q$ that acts elliptically on $T$. Since the cover $\tilde S\to S$ is characteristic, any homeomorphism of $S$ lifts to a homeomorphism of $\tilde S$, and hence the image of $\gamma$ under any homeomorphism of $S$ also is at most one-ended, and therefore acts elliptically on $T$. The point is then that this gives a \emph{lot} of elements of $Q$ that act elliptically, enough so that Serre's criterion (see \ref{sec: trees}) implies that the action of $Q$ on $T$ has a global fixed point, a contradiction.

There are some complications in extending this argument to the setting of noncompact finite-type surfaces. For instance, one needs a version of Stallings's Theorem that works relative to the peripheral groups, in a certain sense. Luckily, there is indeed a relative version of Stallings's theorem due to Swarup \cite{swarup1977relative}, but it is phrased in terms of group cohomology. In order to apply Swarup's theorem in our context, we wrote an Appendix (see \S \ref{sec: appendix}) translating between the group cohomological interpretation of ends and the version that is natural for our paper. We could not find a readable treatment of this in the existing literature, and we hope this can be a useful reference for others. 

    \medskip

Theorem \ref{thm: not locally finite}, in which we identify non-locally finite characteristic covers of infinite-type surfaces, is proved by exhausting by finite-type surfaces and using Theorem \ref{thm: finite type}, in a way analogous to the infinite-type case of Proposition \ref{Prop:UAC}.

\subsection{Outline of the Paper}

In \S \ref{background} we review the necessary background. We define the surfaces mentioned in the introduction, set conventions, and review actions on trees, ends, Stallings's Theorem, and the classification of infinite-type surfaces. In \S \ref{sec: char covers} we study characteristic covers of surfaces, proving Theorem \ref{thm: finite type}, Theorem \ref{thm: not locally finite}, Proposition \ref{Prop:UAC} and Theorem \ref{thm: finite homology covers} in subsections \ref{sec: finite type}, \ref{sec:pfchar}, \ref{sec:UAC} and \ref{thm: finite homology coverssec}, respectively. In \S \ref{sec: arbitrary covers} we prove Proposition \ref{Prop:arbitrary covers}, and then \S \ref{sec: appendix} is our appendix on ends via group cohomology.

\subsection{Acknowledgements}
We thank Sumanta Das for point out Goldman's work and S\'ebastien Alvarez for several other references related to Proposition \ref{Prop:arbitrary covers}. We also thank the anonymous referee for their comments, which have greatly improved the accuracy and readability of the paper.

IB was partially supported by NSF CAREER award 1654114. 
TC was partially supported by MSCA grant 101107744–DefHyp.
NGV was partially supported by NSF award 2212922 and PSC-CUNY awards 66435-00 54 and 67380-00 55. JT was partially supported by NSF award 2304920.

\section{Background}
\label{background}
\subsection{Conventions}
All surfaces within this work are assumed to be connected, second countable, orientable, and possibly with boundary. We say the surface is \emph{borderless} if it has no boundary, and remind the reader that some of the surfaces which appear in our proofs may have non-empty boundary with possibly non-compact boundary components. All borderless surfaces, whether finite or infinite type, are classified by their genus and the topology of their space of ends, see \Cref{thm: classification}. We let $S^b_{g,p}$ denote the finite-type surface with $b$ boundary components, genus $g$, and $p$ punctures.

There are a handful of special infinite-type surfaces that occur in our main theorem statements:

\begin{enumerate}
\item The \textit{flute surface} is homeomorphic to $\BR^2 \smallsetminus \BZ^2$. 
\begin{figure}[h]
    \centering
       \includegraphics[width=.8\textwidth]{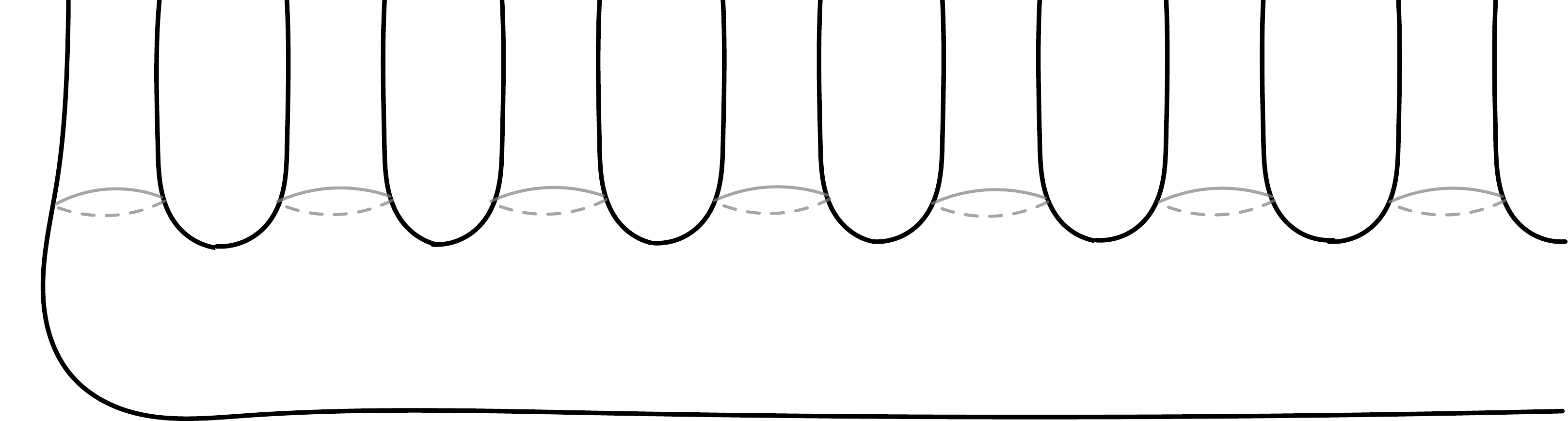}
\end{figure}

\item The \textit{Loch Ness monster surface} is the unique borderless orientable infinite-type surface with one end.
\begin{figure}[h]
  \centering
    \includegraphics[width=.8\textwidth]{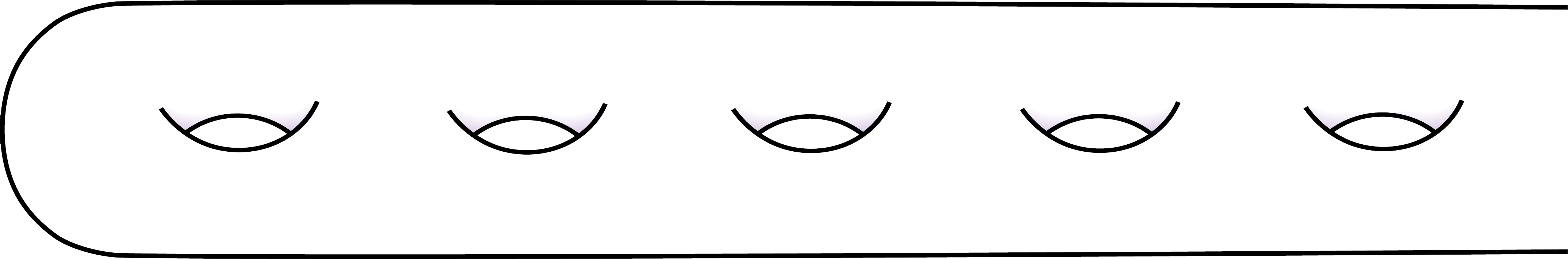}
\end{figure}

\item The \textit{spotted Loch Ness monster surface} is obtained from the Loch Ness monster surface by removing an infinite, discrete, closed set of points.
\begin{figure}[h]
    \centering
     \includegraphics[width=.8\textwidth]{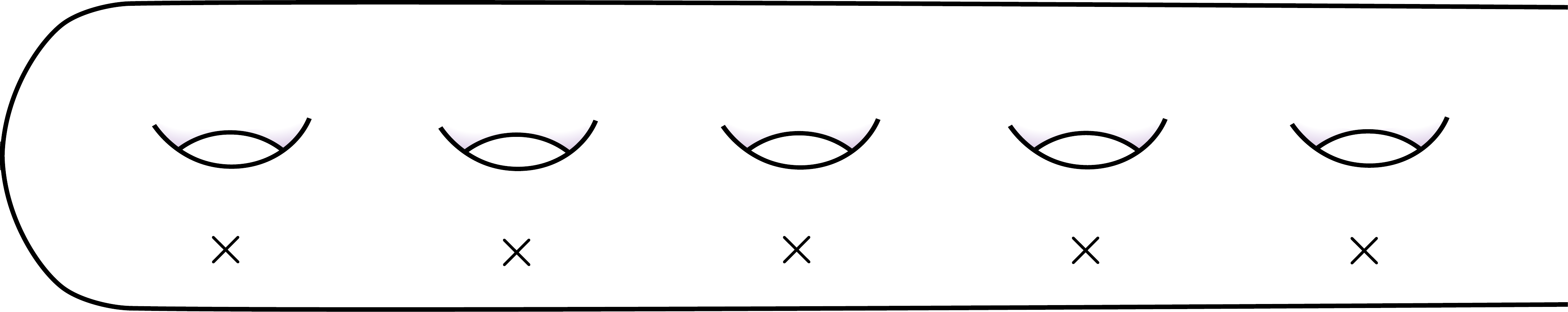}
     \end{figure}
     
\item The \textit{Cantor tree surface} is homeomorphic to $S^2 \smallsetminus \mathcal{C}$, where $\mathcal{C}$ is an embedding of a Cantor set. 
\item The \textit{blooming Cantor tree surface} is the surface with infinitely many genus whose ends space, $\CE(S)$, and ends space accumulated by genus, $\CE_g(S)$, satisfy $\CE_g(S) = \CE(S) \cong \mathcal{C}$. 

\begin{figure}[h]
     \begin{minipage}{\textwidth}
     \centering
        \includegraphics[width=.4\textwidth]{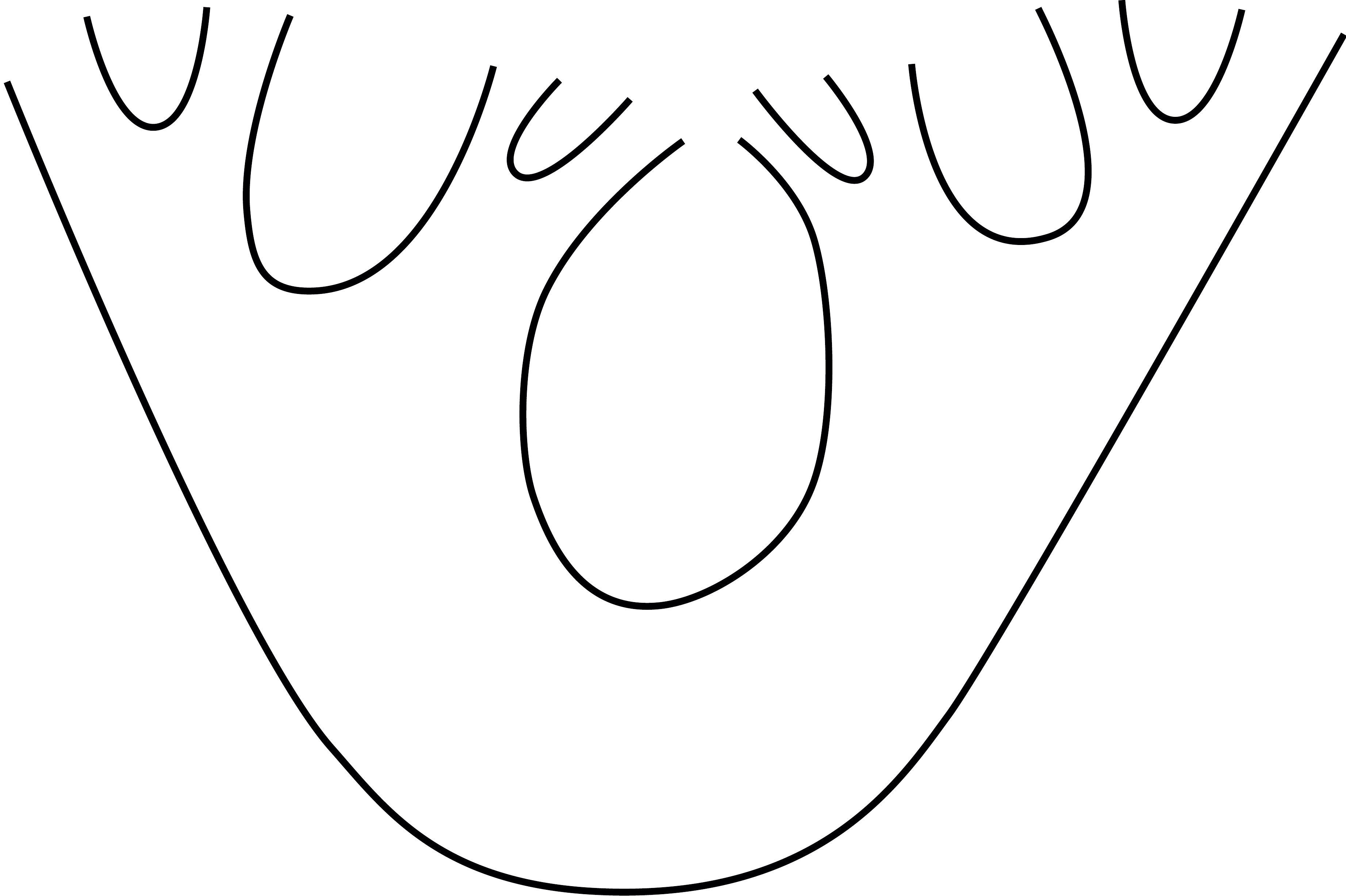} 
         \hspace{1em}\includegraphics[width=.4\textwidth]{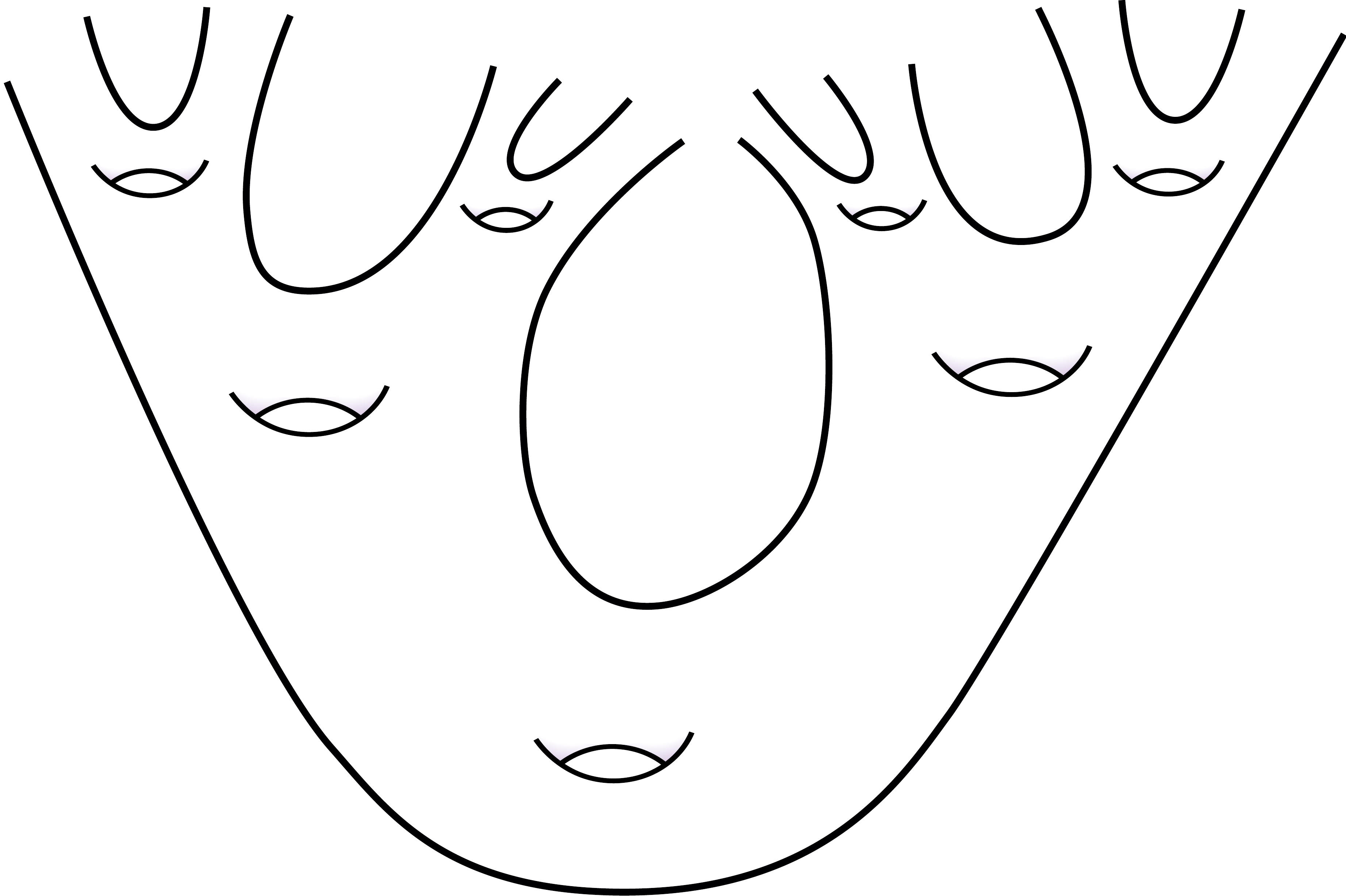}\end{minipage}
\end{figure}

\end{enumerate}

\subsection{Ends of a topological space}
\label{sec: ends}
 Given a topological space $X$, its \textit{space of ends} is defined to be \[ \CE(X) = \lim_{\longleftarrow} \pi_0( X \smallsetminus K),\] where the inverse limit is taken over all compact subsets $K \subset X$. 
 
 For a finitely generated group $G$, consider its Cayley graph $Cay(G,\CS)$ with respect to a finite generating set $\CS$. We define the space of ends of $G$ by $\CE(G) = \CE({Cay}(G,\CS))$. The following proposition shows that $\CE(G)$ is independent of the choice of a finite generating set of $G$.

\begin{proposition}{\cite[Prop. 8.29]{BridsonHaefliger99}}\label{prop:qi-ends}
    Given two proper geodesic metric spaces $X_1$ and $X_2$, then any quasi-isometry $f: X_1 \to X_2$ induces a homeomorphism $f_\CE: \CE(X_1)\to \CE(X_2)$. 
   
\end{proposition}

 \subsection{Classification of surfaces}\label{sec:classification}
Let $S$ be an orientable surface with genus $g(S)\in \{0,\ldots,\infty\}$ and end space $\CE(S)$. Recall that $\CE(S)$ is always compact, totally disconnected, and metrizable. 
We may take a neighborhood of an end to be an open subset $U$ of $S$ containing the end. An end is \textit{planar} if it has a neighborhood that embeds into the plane, and \textit{accumulated by genus} otherwise. Let $\CE_g(S) \subset \CE(S)$ be the closed subset consisting of all ends that are accumulated by genus. The following classification theorem was first proved by K\'erekjarto, and the proof was later simplified  by Richards, see \cite{richards1963classification}.

  \begin{theorem}\label{thm: classification}
    Two borderless, orientable surfaces $S$ and $T$ are homeomorphic if and only if $g(S)=g(T)$ and there is a homeomorphism of pairs $$(\CE(S),\CE_g(S)) \to (\CE(T),\CE_g(T)).$$
    Moreover, for every $n\in \{0,\ldots,\infty\}$ and for every compact, totally disconnected metrizable space $X$, and closed subset $C\subset X$, there is an orientable surface $S$ such that $g(S)=n$ and $(\CE(S),\CE_g(S))$ is homeomorphic to $(X,C)$.
  \end{theorem}

Below, we'll call the pair $(\CE(S),\CE_g(S))$ the \emph{genus-marked end space} of $S$.

\subsection{Preliminaries on covers}
By a \textit{cover} of $S$, we mean a space $\tilde S$ equipped with a surjective map $\pi: \tilde S \to S$ so that every point $p\in S$ has a neighborhood $U$ so that for each component $U_\alpha$ of $\pi^{-1}(U)$, the restriction $\pi|_{U_\alpha}: U_\alpha \to U$ is a homeomorphism. We will occasionally use the word ``cover" just to refer to the covering surface $\tilde S$. The Galois correspondence associates subgroups of $\pi_1(S)$ with (equivalence classes of) connected covers of $S$. 

\begin{remark}
    Suppose $S$ is a surface with non-abelian $\pi_1(S)$, and $H\leq \pi_1(S)$ is a nontrivial, infinite-index, normal subgroup. Then the cover $\tilde S$ corresponding to $H$ is an infinite-type surface. 
\end{remark}

We record a general proof that the cover of a surface corresponding to the fundamental group of an essential subsurface is homeomorphic to the interior of the subsurface. For us, a subsurface of a surface $S$ is always a properly embedded, closed submanifold of $S$. Given a $\pi_1$-injective subsurface \( Y \) of a surface \( S \) and a point \( x \in Y \), we identify \( \pi_1(Y,x) \) with its image under the monomorphism \( \pi_1(Y,x) \to \pi_1(S,x) \) induced by the inclusion \( Y \hookrightarrow S \). When necessary, we will keep track of base points, otherwise we omit them and simply write \( \pi_1(Y) \) and view it as a subgroup of \( \pi_1(S) \).

\begin{lemma}\label{lem:coverofsubsurf}
Let $Y$ be a $\pi_1$--injective subsurface of a borderless surface $S$. Then the cover $\tilde Y$ of $S$ corresponding to $\pi_1(Y)$ is homeomorphic to the interior of $Y$. 
\end{lemma} 
\begin{proof} 

Let $\pi: \tilde Y \to S$ be the cover corresponding to $\pi_1(Y)$ and observe that the inclusion $Y \to S$ lifts to an embedding  $Y \to \tilde Y$, which induces an isomorphism $\pi_1(Y) \to \pi_1(\tilde Y)$. Therefore, any loop in $\tilde{Y}$ is a homotopic to one in the interior of $Y$. This implies that each boundary component of $Y$ is separating---otherwise, there is an essential loop of $\tilde Y$ which essentially intersects a component of $\partial Y$. Further, each component of $\tilde{Y} \smallsetminus Y$ is topologically a product. This shows $\tilde Y\cong  Y\cup (\partial Y\times [0,\infty))$. \qedhere
\end{proof}

\subsection{Group actions on trees}
\label{sec: trees}
Let $T$ be a simplicial tree, and let $G$ be a group acting on $T$. By convention, we assume all such group actions are by simplicial isomorphisms and without edge inversions. Here, an \emph{edge inversion} is a tree automorphism that leaves invariant some edge but reverses its orientation.  An element $g\in G$ acts \emph{elliptically} on $T$ if $g$ fixes a vertex of $T$, and \emph{hyperbolically} otherwise, in which case it translates along a geodesic axis. 
We direct the reader to Serre's book on Trees \cite{serre2002trees} for more information about group actions on trees.

\subsubsection{Stallings's Theorem and its relative version}
  
\label{sec: stallings}

  \begin{theorem}[Stallings's Theorem, \cite{stallings1971group}]

    Suppose $Q$ is a finitely generated group with more than one end. Then $Q$ admits an action on a simplicial tree $T$ with finite edge stablizers, no edge inversions, and no global fixed point.
    
  \end{theorem}

We also need the following relative version of Stallings's Theorem. 

\begin{theorem}[Relative Stallings's Theorem]
Let $Q$ be a finitely generated group and let $H_1,\ldots,H_m$ be finitely generated subgroups of $Q$. Suppose there is a nonconstant, continuous function $f: \CE(Q) \to \BZ$ that is constant on the boundary $\overline {H_ig} \cap \CE(Q)$ of each right coset $H_i g$, where $i=1,\ldots,m$ and $g\in Q$. Then $Q$ admits an action on a simplicial tree $T$ with finite edge stabilizers, no inversions, no global fixed point, and where for each $i$, the action of $H_i$ on $T$ has a global fixed point.\label{Thm:relativestallings}
\end{theorem}

Theorem~\ref{Thm:relativestallings} is really a result of Swarup \cite{swarup1977relative}, but Swarup's result is stated in terms of group cohomology with group ring coefficients. In Appendix~\ref{sec: appendix}, we give a gentle introduction to ends via group cohomology and show how the formulation above follows from Swarup's result. 

The intuition behind the relative version is as follows. Let ${\sim}$ be the smallest equivalence relation on $\CE(Q)$ such that 
\begin{enumerate}
\item for every right coset $H_j g$, all points in $\overline {H_j g}\cap \CE(Q)$ are equivalent,
\item ${\sim}$--equivalence classes are closed in $\CE(G)$.
\end{enumerate}
Taking the quotient $\CE(Q)/{{\sim}}$, we are continuously collapsing the ends of all these cosets. The hypothesis in Theorem \ref{Thm:relativestallings} is that $\CE(G)/{\sim}$ has more than one element, which intuitively means that $Q$ has more than one end even up to collapsing cosets of the $H_i$. So, there should be an action of $Q$ on a tree $T$, as in Stallings's Theorem, in which the $H_i$ are not really used.

In this paper, we use Theorem~\ref{Thm:relativestallings} in the following rather explicit context. The group acting will be the deck group $Q$ of a regular cover $\tilde S \to S$, where $S$ is a compact surface with boundary satifying $\tilde \CE(Q) \cong \tilde \CE(\tilde S)$. The subgroups $H_i<Q$ will be stabilizers of noncompact boundary components of $\tilde S$. Deleting these boundary components creates a new surface $\tilde S_\circ$, in which the two ends of each such boundary component become a single end of $\tilde S_\circ$. So, the end space $\CE(\tilde S_\circ) $ is the same as the quotient $\CE(Q)/{\sim}$ discussed in the previous paragraph, and the hypothesis in Theorem~\ref{Thm:relativestallings} is that $\tilde S_\circ$ has more than one end. In the paper, we phrase all of this in a way that appeals directly to the statement of Theorem~\ref{Thm:relativestallings}, though, rather than talking about the relation ${\sim}$ above.

In addition to Stalling's Theorem and our formulation of its relative version, our main arguments in Section \ref{sec: finite type} uses Serre's criterion to cook up a proof by contradiction to ensure the surface $\tilde S_\circ$ is one-ended. 
  \begin{lemma}[Serre's Criterion, see Corollary~2 in \S 6.5 of \cite{serre2002trees}]\label{Lem:Serre} 

    Let \( G \) be a group acting on a tree \( T \). If \( G \) admits a
    generating set \( \{a_1, \ldots, a_n\} \) such that \( a_i \) and \(
    a_ia_j \) act elliptically on \( T \) for each distinct \( i \) and \(
    j \), then the $G$-action has a global fixed point.
  \end{lemma}

\section{Characteristic covers}
\label{sec: char covers}

In this section we study characteristic covers of surfaces, proving Theorems \ref{thm: finite type}, \ref{thm: not locally finite}, Proposition~\ref{Prop:UAC}, and Theorem~\ref{thm: finite homology covers} in sections \ref{sec: finite type}, \ref{sec:pfchar}, \ref{sec:UAC}, and \ref{thm: finite homology coverssec}, respectively.

\subsection{Finite-type surfaces}
\label{sec: finite type}

Here, we prove the following theorem from the introduction.

\begin{named}{Theorem \ref{thm: finite type}}
    Any infinite-degree, geometrically characteristic cover of an orientable, borderless, finite-type surface is either a disc, the flute surface, the Loch Ness monster surface, or the spotted Loch Ness monster surface.
\end{named}

Throughout the section, we use the following notation: 
Given a surface \( S \), we let \( S_\circ \) denote the surface obtained by removing every non-compact boundary component of \( S \). We'll prove the following statement.

\begin{proposition}\label{prop: rephrase finite type}
Let \( S \) be a compact, orientable surface.
If $\tilde S \to S$ is a geometrically characteristic cover of infinite degree, then  $\tilde S_\circ$ has one end. 
\end{proposition}

Assuming Proposition \ref{prop: rephrase finite type}, here is how to derive Theorem \ref{thm: finite type}.

\begin{proof}[Proof of Theorem~\ref{thm: finite type}] 

Note that any finite-type surface can be identified as the interior of some compact surface $S$. 
Further, a geometrically characteristic cover of the interior of $S$ is the interior of a geometrically characteristic cover $\tilde{S} \to S$. 
Thus, the goal is to show that the interior of $\tilde{S}$ is one of the four surfaces listed. 

Let $Q$ be the associated deck group of $\tilde{S}$. 
First assume that \( \tilde S \) has no compact boundary components, so $\int(\tilde S) = \tilde S_\circ $.
As $Q$ is infinite and  acts cocompactly on $\tilde S$, the genus of $\int(\tilde S)$ is either $0$ or infinity, so the classification of surfaces (Theorem \ref{thm: classification}) implies that it is either the disc or the Loch Ness monster surface. 
Now, suppose that $\tilde S$ has a compact boundary component.
As $Q$ is infinite, $\tilde{S}$ has infinitely many compact boundary components. 
Gluing discs to each compact boundary component, we obtain a surface without boundary, call it $Y$.
By Proposition~\ref{prop: rephrase finite type}, \( Y \) has a single end, and again the classification of surfaces (Theorem \ref{thm: classification}) implies that $Y$ is either the disc or the Loch Ness monster surface. 
To finish, observe that $\int(\tilde S)$ is obtained from $Y$ by puncturing all the infinitely many discs, so then $\int(\tilde S)$ is either the flute surface or the spotted Loch Ness monster surface.
\end{proof}

In the proof of Proposition \ref{prop: rephrase finite type}, we will start by assuming $\tilde S_\circ$ has more than one end. We then use a relative version of Stallings' Theorem to produce an action on a tree $T$, and decompose $\tilde S$ equivariantly into subsurfaces essentially corresponding to the vertices of $T$. We use the induced decomposition of $S$ to construct elements of $\pi_1 (S)$ that act elliptically on $T$, then use the fact that the cover is characteristic to construct many more elements that act elliptically, in fact so many that the action on $T$ must have a global fixed point, a contradiction. 

For clarity of exposition, though, before starting the proof properly in \S \ref{sec:proofofprop}, we isolate some of its components into a series of preliminary results. First, the following will allow us to construct the action on the tree referenced above.

\begin{lemma}[A peripherally elliptic action]\label{relativeends}
Let $\pi: \tilde S \to S$ be a regular cover of a compact surface \( S \).
If $\tilde S_\circ$ has more than one end, then there is an action of the deck group $Q$ on a tree $T$ with no edge inversions and no global fixed point, such that the projection to $Q$ of every peripheral element of $\pi_1(S)$ acts elliptically on $T$.
\end{lemma}

\begin{proof}
Let $\bar \iota : \CE(\tilde S) \to \CE(\tilde S_\circ)$ be defined as follows: given $\xi \in \CE(\tilde S)$, take a sequence of points $x_i \in \tilde S_0\subset \tilde S$ with $x_i \to \xi$, and define  $\bar \iota(\xi)$ to be the limit in $\CE(\tilde S_\circ)$. 
It is readily verified that \( \bar \iota \) is a well-defined continuous surjection.

As the action of \( Q \) on \( \tilde S \) is co-compact and proper, the end space of \( Q \), denoted \( \CE(Q) \), is homeomorphic to \( \CE(\tilde S) \).
If $\tilde S_\circ$ has more than one end, then there is a non-constant continuous function $\phi : \CE(\tilde S_\circ) \to \BZ$. 
Consider the composition $$\CE(Q) \overset{\cong}{\to} \CE(\tilde S) \overset{\bar \iota}{\to} \CE(\tilde S_\circ) \overset{\phi}{\to} \BZ.$$
Pick cyclic subgroups of $\pi_1 (S)$ representing the boundary components, and let $H_j$ denote their projections to $ Q $.
For any right coset $H_j g$, the intersection $\overline{H_j g}\cap \CE(Q)$ is either empty (if $H_j$ is finite, which happens when that boundary component lifts to a closed curve in $\tilde S$) or projects to the endpoints in $\CE(\tilde S)$ of some non-compact boundary component of $\tilde S$, both of which (if they are distinct) map to the same element of $\CE(\tilde S_\circ)$. 

We have established that $\phi$ is non-constant but that it is constant on $\overline{H_j g}\cap \CE(Q)$ for each $j$ and each right coset $H_jg$. 
Therefore,  by Theorem~\ref{Thm:relativestallings} there exists an action of $Q$ on a tree $T$, with no edge inversions and finite edge stabilizers, and no global fixed point, such that all the subgroups $H_j$ act elliptically. 
Every peripheral element of $\pi_1 (S)$ projects to an element of $Q$ that is conjugate into one of the $H_j$'s and hence acts elliptically on $T$ as desired.
 \end{proof}

Given a peripherally elliptic action of $Q$ on a tree $T$, the following proposition constructs a $Q$-equivariant collection of simple closed curves in $\tilde S$ corresponding (perhaps multiple-to-one) to the edges of $T$. Cutting along these curves decomposes $\tilde S$ into subsurfaces corresponding to the vertices of $T$.
 
 \begin{proposition}
 \label{prop: equiv map}
Let $\pi: \tilde S \to S$ be a regular cover of a compact surface \( S \)  with deck group \( Q \).
Suppose \( Q \) acts on a tree \( T \) with finite edge stabilizers, no edge inversions, no global fixed point, and such that the projection to \( Q \) of every peripheral element of \( \pi_1(S) \) acts elliptically on \( T \). 
Then there exists a $Q$-equivariant map $\tilde f : \tilde S \to T$ such that  each component of \( \partial \tilde S \) maps to a vertex under \( \tilde f \).
Moreover, if \( Z \) is the set of edge midpoints in \( T \), then the preimage \( \tilde C := \tilde f^{-1}(Z) \) is a properly embedded $1$-submanifold of $\tilde S$, each component of which is a simple closed curve, and its projection \( C := \pi(\tilde C) \) is a compact $1$-submanifold of $S$ each component of which is non-peripheral and essential. 
\end{proposition}

Above, there may be components of $C$ that are parallel, which is why we call $C$ a $1$-submanifold and not a multicurve.

\begin{proof}
We begin by constructing a prototype for the map \( \tilde f \), which we call \( \tilde f_1 \). 
Following Shalen \cite{shalen2002representations}, we construct the $Q$-equivariant map $\tilde f_1\co \tilde S \to T$ as follows. 
Fix a triangulation $\tau$ of $S$, and consider its associated $Q$-invariant triangulation $\tilde \tau$ of $\tilde S$.
Let $\tilde X_0 \subset T^{(0)}$ contain  a single representative for each orbit equivalence class of $\tilde \tau^{(0)}$, and pick a map $$\tilde h^{(0)}\co \tilde X_0 \to T^{(0)}$$ that is arbitrary, except that for every component $\beta\subset \partial \tilde S$, we require that $\tilde h^{(0)}$ takes $\tilde X_0 \cap \beta$ to a single vertex of $T$ that is fixed by the stabilizer $Q_\beta \subset Q$; this can be done as $Q_\beta$ is the projection to $Q$ of a peripheral $\BZ$-subgroup of $\pi_1(S)$, so acts with a global fixed point on $T$ by assumption. (Note that $Q_\beta$ must fix a \emph{vertex}, since it acts simplicially without edge inversions.)

\begin{figure}[t]
\centering
\includegraphics{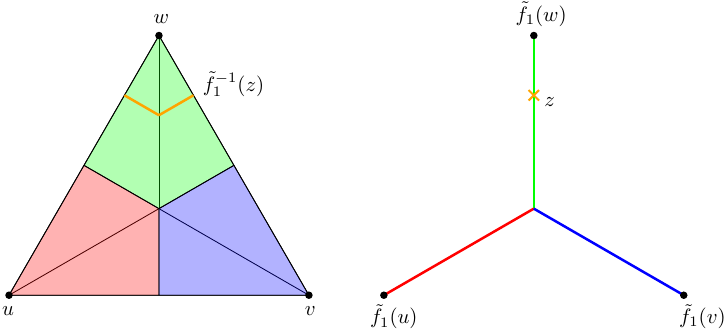}
\caption{The map \( \tilde f_1 \) from Proposition~\ref{prop: equiv map} restricted to a triangle in the triangulation \( \tilde \tau \); the green, blue, and red regions map onto the green, blue, and red segments, respectively, in the tripod; the restriction to any one of the six subtriangles is an orthogonal projection map.
The preimage of an edge midpoint \( z \in T \) is shown in orange.}
\label{fig: transverse}
\end{figure}

There is a unique $Q$-equivariant extension of $\tilde h^{(0)}$ to a map $$\tilde f_1^{(0)}\co \tau^{(0)} \to T^{(0)}.$$ 
The map \( \tilde f_1^{(0)} \) is readily extended to a piecewise linear map \( \tilde f_1 \co \tilde S \to T \) (see Figure~\ref{fig: transverse}). 
Observe that, under this piecewise linear extension, each boundary component $\beta$ of $\partial \tilde S$ is mapped to a point by $\tilde f_1$, as its vertices are all mapped to the same point by construction. 

Let \( Z \) be the set of all edge midpoints of \( T \), let \( \tilde C_1 = \tilde f_1^{-1}(Z) \), and let \( C_1 = \pi(C_1) \).
We claim that each component \( \tilde C_1 \) (resp., \( C_1 \)) is a simple closed curve. 
By the construction of \( \tilde f_1 \), we know that the components of \( \tilde C_1 \) form a locally finite family of pairwise-disjoint 1-submanifolds of \( \tilde S \). 
As \( \pi \) is a local homeomorphism and \( \tilde C_1 \) is invariant under the action of \( Q \), we have that \( C_1 \) is a compact 1-submanifold of \( S \); in particular, each component of \( C_1 \) is a simple closed curve. 
Now, let \( \tilde c \) be a component of \( \tilde C_1 \).
Then \( \pi(\tilde c) \) is compact, and as the action of \( Q \) on \( T \) has finite edge stabilizers, the stabilizer of \( \tilde c \) is finite, implying that \( \tilde c \) is compact. 
Therefore, \( \tilde C_1 \) is a union of simple closed curves. 
We now need to edit \( \tilde f_1 \) so that no component of \( \tilde C_1 \) nor \( C_1 \)  is peripheral or bounds a disk. 

First, suppose that there is a component of \( C_1 \) that bounds a disk. 
Let \( \Delta_1, \ldots, \Delta_k \) be the outermost disks bounded by components of \( C_1 \), so that if \( c \) is a component of \( C_1 \) bounding a disk, then there exists a unique \( j \in \{1, \ldots, k \} \) such that \( c \subset \Delta_j \). 
Choose a lift \( \tilde \Delta_j \) in \( \tilde S \) for each of the \( \Delta_j \), and for \( q \in Q \), let \( \tilde \Delta_j^q = q \cdot \tilde \Delta_j \). 
Let $e_j^q$ be the edge with vertices $v_j^q,w_j^q \in T$ whose midpoint is $\tilde{f_1}(\tilde \partial \Delta_j^q)$.  
By the construction of $C_1$, all the $1$-cells in $\tilde \tau$ that intersect $\tilde \partial \Delta_j^q$ are mapped to segments in $T$ that all contain the edge $e_j^q$.  
If the $0$-cells in \( \tilde \Delta_j^q \) are mapped into the $v_j^q$ side of $T\smallsetminus\{e_j^q\}$, we can define a new map $$\tilde f_2^{(0)} : \tilde \tau^{(0)} \to T ^{(0)}$$ that maps the 0-cells in \( \Delta_j^q \) to $w_j^q$, but agrees with \( \tilde f_1^{(0)} \) otherwise.
We can now extend \( \tilde f_2^{(0)} \) to a \( Q \)-invariant piecewise linear map \( \tilde f_2 \co \tilde S \to T \) in such a way that \( \tilde f_2^{-1}(z) = \tilde f_1^{-1}(z) \) for any \( z \in Z \smallsetminus (\bigcup \tilde \Delta_j^q) \). 
Setting \( \tilde C_2 = \tilde f_2^{-1}(Z) \) and \( C_2 = \pi( \tilde C_2 ) \), we have that \( C_2 \subset C_1 \) and no component of \( C_2 \) bounds a disk. 

Now, suppose that \( C_2 \) contains a peripheral component. 
Then each peripheral component of \( C_2 \) bounds an annulus with a component of \( \partial S \). 
Let \( A_1, \ldots, A_m \) be the outermost annuli, so that if \( c \) is a component of \( C_2 \) bounding an annulus with a component of \( \partial S \), then there exists a unique \( j \in \{1, \ldots, m\} \) such that \( c \subset A_j \). 
Each lift of an \( A_j \) is again an annulus in \( \tilde S \) as the component of \( \tilde C_2 \) are compact. 
We now proceed identically as we did above with the disks \( \Delta_j \) replaced with the annuli \( A_j \). 
The result is a \( Q \)-equivariant map \( \tilde f \co \tilde S \to T \) so that \( \tilde C := \tilde f^{-1}(Z) \) and \( C := \pi(\tilde C) \) have neither peripheral nor inessential components. 
 \end{proof}

Combining Lemma~\ref{relativeends} and  Proposition~\ref{prop: equiv map}, we obtain:

\begin{corollary} \label{cor:curve}
Let $\pi: \tilde S \to S$ be a regular cover of a compact surface.
If $\tilde S_\circ$ has more than one end, then there is a non-peripheral, essential, simple closed curve $\gamma$ on $S$ such that every component of $\pi^{-1}(\gamma)$ is compact.

\end{corollary}
\begin{proof}
    Since the action of the deck group $Q$ produced by Lemma~\ref{relativeends} has no global fixed point, the image of the map in Proposition~\ref{prop: equiv map} has to contain an edge midpoint, the $1$-submanifold $C\subset S$ it produces is non-empty. Take $\gamma$ to be any component of $C$. Any component of $\pi^{-1}(\gamma)$ is a component of $\tilde C$, and hence is compact.
\end{proof}

Given a properly embedded $1$-submanifold \( C \) on a surface \( S \), we write \( S \dsm C \) to refer to the surface obtained by cutting \( S \) along \( C \), that is, the surface obtained by removing an open regular neighborhood of \( C \) from \( S \). 

\begin{lemma}[Cutting covers into one-ended pieces]
\label{lem:multicurve}
If $\pi: \tilde{S} \to S$ is a regular cover of a compact surface and $C\subset S$ is a compact $1$-submanifold such that every component of $\pi^{-1}(C) \subset \tilde S$ is compact.
 Then there exists a compact $1$-submanifold $D\subset S$ containing $C$ such that
\begin{itemize}
    \item every component of $\pi^{-1}(D)$ is compact, and
    \item every component of $\tilde S_\circ \dsm \pi^{-1}(D)$ has at most one end.
\end{itemize} 
\end{lemma}

\begin{proof}
It suffices to show that if some component of $\tilde S_\circ \dsm \pi^{-1}(C)$ has more than one end, then there is an essential, nonperipheral simple closed curve $\gamma  \subset S \dsm C$ such that every component of $\pi^{-1}(\gamma) \subset \tilde S$ is compact. For then we can add such curves recursively to $C$ to produce $D$ as required. Since the number of isotopy classes in a collection of pairwise disjoint curves on $S$ is bounded, this process terminates.

So, let \( \Sigma \) be a component of $S \dsm C$, and assume that we have a component $\tilde \Sigma\subset \tilde S \dsm \tilde C$ such that $\tilde \Sigma_\circ$ has more than one end. We claim that the restriction of \( \pi \) to \( \tilde \Sigma \) is a regular covering map $ \tilde \Sigma \to \Sigma.$ 
 Indeed, if two points $p,q \in \tilde \Sigma$ have the same projection to $\Sigma$, then there is a deck transformation $g$ of $\tilde S \to S$ with $g(p)=q$, and since $\tilde C$ is invariant under the deck group, we have $g(\tilde \Sigma)=\tilde \Sigma$, implying $g$ restricts to a deck transformation of $\tilde \Sigma \to \Sigma$ taking $p$ to $q$.  
 
 Corollary~\ref{cor:curve} then says that there is some essential, non-peripheral curve in $\gamma \subset \Sigma$ such that every component of $\pi|_{\tilde \Sigma}^{-1}(\gamma)$ is compact. But since the cover $\pi : \tilde S \to S$ is regular, this implies that every component of $\pi^{-1}(\gamma) \subset \tilde S$ is compact, so we are done.
\end{proof}

\subsubsection{Proof of Proposition~\ref{prop: rephrase finite type}}\label{sec:proofofprop}
We argue by way of contradiction: assume that \( \tilde S_\circ \) has more than one end. 
Let \( Q \) be the deck group associated to the covering \( \pi \). 
Apply Proposition~\ref{prop: equiv map} to obtain a tree \( T \), a peripherally elliptic action of \( Q \) on \( T \) with no global fixed point, and the 1-submanifolds \( \tilde C \) and \( C  \) of \( \tilde S \) and \( S \), respectively, and apply Lemma~\ref{lem:multicurve} to get the $1$-submanifold \( D \) on \( S \).
Recall that the action of \( Q \) on \( T \) has finite edge stabilizers, no edge inversions, no global fixed point, and every element of \( Q \) that is the image of a peripheral element in \( \pi_1(S) \) acts elliptically on \( T \).
The contradiction will arise by showing that the action of \( Q \) has a global fixed point. 

We say that a closed curve $\gamma$ is \emph{at most one ended in $\tilde S_\circ$} if for some (equivalently, any) component $\tilde \gamma$ of $ \pi^{-1}(\gamma)$, we have that $\tilde \gamma$ is either compact or it is non-compact and its two ends go out the same end of $\tilde S_\circ$. 
In particular, any closed curve  disjoint from \( D \)  is at most one ended in $\tilde S_\circ$.  

We say two closed curves \( \gamma_1 \) and \( \gamma_2 \) in \( S \) have the same homeomorphism type if there exists a homeomorphism \( h \co S \to S \) such \( h(\gamma_1) = \gamma_2 \). 
As \( \pi \) is geometrically characteristic, every homeomorphism \( h \co S \to S \) lifts to a homeomorphism \( h \co \tilde S \to \tilde S \).
Also note that the restriction of any homeomorphism \( \tilde S \to \tilde S \) to \( \tilde S_\circ \) induces a homeomorphism \( \tilde S_\circ \to \tilde S_\circ \). 
Therefore, if \( \gamma \) is a closed curve that has the same homeomorphism type as a closed curve  disjoint from \( D \),  then \( \gamma \) is at most one ended in \( \tilde S_\circ \).

\medskip

Now, we claim that if \( \gamma \) is an essential closed curve on \( S \) that is at most one ended, then \emph{$\gamma$ acts elliptically on $T$}, i.e., any element of \( \pi_1(S) \) in the free homotopy class of \( \gamma \) acts elliptically on \( T \), where we let \( \pi_1(S) \) act on \( T \) via the action factoring through \( Q \). 

To see this, pick some $g\in Q$ that is the projection of an element of $\pi_1(S)$ in the free homotopy class of $\gamma$, and hoping for a contradiction, \emph{suppose that $g$ acts hyperbolically on $T$}. Now $g$ stabilizes some component $\tilde \gamma \subset \pi^{-1}(\gamma)$. Since $g$ acts hyperbolically on $T$, the image $\tilde f(\tilde \gamma)$ is a proper bi-infinite path in $T$, whose ends go out distinct ends of $T$. 
Pick some edge midpoint $m \in T$ that separates those two ends, and let $M = \tilde f^{-1}(m). $ 
Then $M$ is a union of components of $\tilde C$. 

We claim that  $M$ has finitely many components. If not, then since the projection $C$ has only finitely many components, there are infinitely many elements of $Q$ that send a component of $M$ to a component of $M$, and hence $m \in T$ has infinite stabilizer, contradicting our assumption that $Q$ acts with finite edge stabilizers. 
Hence, $M$ has finitely many components.

The two ends of the path $\tilde \gamma \subset \tilde S$ lie in different components of $\tilde S \smallsetminus M$, since they map to different components of $T\smallsetminus m$. Since $M$ is compact, this means that the two ends of $\tilde \gamma$ converge to distinct ends of $\tilde S_\circ$, \emph{a contradiction.}
This establishes the claim that \( \gamma \) acts elliptically on \( T \). 

\medskip

Thus far, we have established that any essential closed curve on $S$ that has the same homeomorphism type as a curve  disjoint from  $D$ acts elliptically on $T$. 
In the rest of the proof, using this plethora of elliptic elements, we apply Serre's Criterion (Lemma \ref{Lem:Serre}) to show that the action $Q \actson T$ has a global fixed point.

\medskip

First, let us assume that there is no non-separating simple closed curve on $S$ that is disjoint from $D$. 
Then the components of $S\smallsetminus D$ are all planar and are connected together in the pattern of a tree, since if there is a cycle in the component graph, you can realize it as a simple closed curve on $S$ intersecting some component of $D$ once, and then that component of $D$ is non-separating, contradicting that there are no non-separating curves disjoint from $D$. 
It follows that $S$ is planar, and hence its fundamental group can be generated by a set of elements representing all but one of the boundary components of $S$. 
By the construction of the action of \( Q \) on \( T \), all such curves act elliptically on $T$. 

A \emph{boundary product} is a simple closed curve on $S$ that bounds a pair of pants with two boundary components of $S$. Note that all boundary products in $S$ differ by homeomorphisms of $S$.
If we choose orientations for the generators above correctly, all products of pairs of distinct generators are homotopic to boundary products.  Looking within a component of $S\smallsetminus D$ that is a leaf of the associated component tree, we can find a {boundary product} disjoint from $D$, so it follows from the above that \emph{all} boundary products act elliptically on $T$, and hence the action has a global fixed point by Serre's Criterion (Lemma \ref{Lem:Serre}), a contradiction.

\medskip

\begin{figure}[t]
\centering
\includegraphics{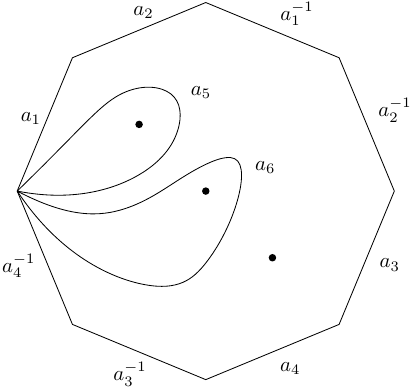}
\caption{Standard generators for \( \pi_1(S) \) for \( S \) on a genus two surface with three punctures are given by the elements \( a_1, \ldots, a_6 \).}
\label{fig:standard_generators}
\end{figure}

Now assume there is a non-separating simple closed curve disjoint from $D$. 
Then all non-separating curves on $S$ act elliptically on $T$. 
Consider a `standard' generating set for $\pi_1(S)$ consisting of two non-separating curves intersecting once for each genus, plus generators for all but one peripheral curve (see Figure~\ref{fig:standard_generators}). 
All these generators act elliptically on $T$. 
A product of a pair of generators is freely homotopic to either a boundary product or a non-separating simple closed curve.   
To use Serre's Criterion again to get a contradiction, it suffices to prove that all boundary products on $S$ act elliptically on $T$.

\begin{figure}[h]
\centering
\includegraphics{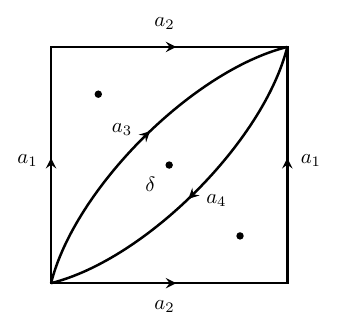}
\caption{A generating set for a thrice-punctured torus consisting of non-separating curves, namely \( a_1, a_2, a_3, a_4 \), where products of pairs of distinct generators are either non-separating, or are freely homotopic to $\delta$. For a similar example on a twice-punctured torus, just delete one of the two boundary components that is not $\delta$.}
\label{fig:torus_generators}
\end{figure}

Fix a boundary product $\gamma$. Recall that this means that $\gamma$ bounds a pair of pants in $S$. Call $\delta$ and $\delta'$ the two peripheral curves of this pair of pants. Since $S$ has a non-separating curve, it has positive genus. It follows that there is an essential genus one subsurface $X \subset S$ containing $\gamma$ that has either two or three boundary components. Indeed, if $S$ has genus one and only two boundary components, then set $X=S$. Otherwise, we can find a separating curve in $S$ that cuts off a subsurface $X$ of genus one containing $\gamma$, and $X$ only has three boundary components. In either case, both $\delta$ and $\delta'$ are peripheral curves of $X$.

Given such an $X$, we can find a generating set for $\pi_1 (X)$ consisting of non-separating curves, such that every product of distinct pairs of generators is either non-separating or is freely homotopic to $\delta$. See Figure \ref{fig:torus_generators}. Since non-separating curves and peripheral curves in $S$ act elliptically on $T$, by Serre's criterion we have that the action $\pi_1 X \actson T$ has a global fixed point. In particular, $\gamma$ acts elliptically on $T$. Since $\gamma$ was an arbitrary boundary product, all boundary products act elliptically on $T$. 

We have now established that the action of \( Q \) on \( T \) has a global fixed point, a contradiction.
This finishes the proof of Proposition \ref{prop: rephrase finite type}.

\subsection{Non-locally finite covers of infinite-type surfaces}\label{sec:pfchar}

Recall that a subgroup \( N \) of a group \( G \) is \emph{locally finite index} if \( N \cap H \) is finite index in \( H \) whenever \( H < G \) is finitely generated; a cover \( \pi \co \tilde S \to S \) is \emph{locally finite} if \( \pi_*(\pi_1 \tilde S) \) is locally finite index in \( \pi_1(S) \). 
In this section we prove:

 \begin{named}{Theorem \ref{thm: not locally finite}}
    Let $S$ be an orientable, infinite-type borderless surface.
    If $\tilde S \to S$ is a characteristic cover that is not locally finite, then $\tilde S $ is homeomorphic to either a disc, the flute surface, the Loch Ness monster surface, or the spotted Loch Ness monster surface.
\end{named}

We first start with some terminology and three lemmas. Recall that a subgroup \( A \) of a group \( G \) is a \emph{free factor} if there exists another subgroup \( B \) of \( G \) such that \( G \) is the internal free product of \( A \) and \( B \). The  three lemmas below tell us the following: given a non-locally finite characteristic cover \( \pi \co \tilde S \to S \) and a ``large'' compact subsurface \( Y \) of \( S \), the restriction of \( \pi \) to a component \( \tilde Y \) of \( \pi^{-1}(Y) \) is an infinite-sheeted characteristic cover \( \tilde Y \to Y \).

\begin{lemma}\label{lem:freefactor}
Let \( S \) be a noncompact surface, and let \( Y \) be a compact subsurface of \( S \). 
If  each component of \( \partial Y \) is a separating curve and each component of \( S \smallsetminus int(Y) \) is noncompact, then
 \( \pi_1(Y) \) is a free factor of \( \pi_1(S) \).  
\end{lemma}

\begin{proof}
Let $\alpha$ be a separating curve in $S$ with complementary components $\Sigma_1$ and $\Sigma_2$. 
We observe that if neither $\Sigma_1$ nor $\Sigma_2$ is precompact, then $\alpha$ can be extended to a free basis for $\pi_1(\Sigma_1)$ and $\pi_1(\Sigma_2)$. By Van Kampen's theorem, $$\pi_1(S) = \pi_1(\Sigma_1) * \pi_1(\Sigma_2).$$ It is not hard to see from this that any free basis for $\pi_1(Y)$ can be extended to a free basis for $\pi_1(S)$. \qedhere

\end{proof}

\begin{lemma}
\label{lem:noncompact_component}
Let \( S \) be a noncompact surface, and let \( \pi \co \tilde S \to S \) be a non-locally finite characteristic cover.
Then there exists \( R > 0 \) such that if  \( Y \subset S \) is an $\pi_1$-injective subsurface whose fundamental group  contains a free factor of \( \pi_1(S) \) of rank \( R \), then each component of \( \pi^{-1}(Y) \) is not compact.
\end{lemma}

\begin{proof}
Let \( N = \pi_*(\pi_1(\tilde S)) \).
As \( N \) fails to be locally finite index in \( \pi_1(S) \), there exists a finitely generated group \( H < \pi_1(S) \) such that \( [H: H\cap N] \) is infinite. 
Using the fact that \( H \) is finitely generated, we can find a free factor \( A \) of \( \pi_1(S) \) containing \( H \); let \( R \) denote the rank of \( A \). 
Let \( Y \subset S \) be a $\pi_1$-injective subsurface so that \( \pi_1(Y) \) contains a free factor \( B \) of \( \pi_1(S) \) of rank \( R \). 
Then there exists an automorphism \( \varphi \) of \( \pi_1(S) \) such that \( \varphi(A) = B \). 
As \( N \) is characteristic, \( \varphi(N) = N \), implying that \[  [\varphi(H): \varphi(H) \cap N] = [\varphi(H): \varphi(H \cap N)] = [H: H\cap N] . \]
In particular, as \( \varphi(H) < \pi_1(Y) \), we have that the index of \( \pi_1(Y) \cap N \) in \( \pi_1(Y) \) is infinite.  
Given a component \( \tilde Y \) of \( \pi^{-1}(Y) \), the restriction of \( \pi \) to \( \tilde Y \) is a covering \( \tilde Y \to Y \) with corresponding subgroup \( \pi_1(Y) \cap N \), and hence \( \tilde Y \) is an infinite-sheeted cover of \( Y \), implying it is noncompact. 
\end{proof}

\begin{lemma}\label{lem:geomcharsubcover}
Let \( S \) be a noncompact surface, and let \( Y \subset S \) be a compact subsurface such that each component of \( \partial Y \) is a separating curve and each component of \( S \smallsetminus int(Y) \) is noncompact.
If \( N < \pi_1(S) \) is characteristic, then \( N \cap \pi_1(Y)  \) is a characteristic subgroup of \(\pi_1(Y) \). 
\end{lemma}

\begin{proof}
Let \( \varphi \in \mathrm{Aut}(\pi_1(Y) ) \).
By Lemma~\ref{lem:freefactor}, \( \pi_1(Y)  \) is a free factor of \( \pi_1(S) \), and so there exists \( B < \pi_1(S) \) such that \( \pi_1(S) \) is isomorphic to \( \pi_1(Y) * B \). 
We can therefore extend \( \varphi \) to an automorphism \( \hat \varphi \) of \( \pi_1(S) \) by declaring \( \hat \varphi(b) = b \) for all \( b \in B \) and \(  \hat \varphi |_{\pi_1(Y)} = \varphi \). 
So,
\begin{align*}
\varphi(N\cap \pi_1 (Y) ) 	&= \hat\varphi(N\cap \pi_1(Y) ) \\
				&= \hat\varphi(N) \cap \hat\varphi(\pi_1(Y) ) \\
				& = N \cap \pi_1(Y).
\end{align*}
Hence, \( N \cap \pi_1(Y)  \) is characteristic in \( \pi_1(Y)  \). 
\end{proof}

We are now ready to prove Theorem \ref{thm: not locally finite}.

\begin{proof}[Proof of Theorem \ref{thm: not locally finite}]
Let \( S \) be an infinite-type surface, and let \( \bar S \) denote the surface obtained by removing an open annular neighborhood of each isolated planar end of \( S \). That is, \( \bar S \) is obtained by replacing the isolated planar ends of \( S \) with compact boundary components. 

Let \( \pi \co \tilde S \to \bar S \) be an infinite-sheeted characteristic cover. 
We will prove that \( \tilde S_\circ \) is one ended (recall that $\tilde S_\circ$ is obtained by deleting the non-compact boundary components of $\tilde S$).
Note that we are actually interested in the topology of $int(\tilde S)$; however, if \( \tilde S_\circ \) is one ended, then the result is a consequence of the classification of surfaces.  

Let \( K \subset \tilde S_\circ  \) be a compact subset.
By possibly enlarging \( K \), we can assume that no component of \( \tilde S _\circ \smallsetminus K \) is precompact. Our goal is to show $\tilde{S}_\circ \smallsetminus K$ is connected. 

Let $R$ be the constant of Lemma~\ref{lem:noncompact_component}. We can exhaust \( \bar S \) by compact subsurfaces whose boundary components are separating and whose complementary components are not precompact. Thus, we may find two nested compact subsurfaces \( Y \subset X \subset \bar{S} \) with the following properties.
\begin{itemize}
\item The interior of $Y$ contains $\pi(K)$. 
\item Each  boundary component of \( X \) is separating, and each component of \( \bar S \smallsetminus int(X) \) is noncompact, 
\item Each component of $X \smallsetminus int(Y)$ has a fundamental group containing a free factor of rank $R$. 
\end{itemize}

Now let $\tilde Y$ be the component of $\pi^{-1}(Y)$ containing $K$, and let $\tilde X  \supset \tilde{Y}$ be the associated component of $\pi^{-1}(X)$. By Lemma~\ref{lem:geomcharsubcover}, the restriction \( \pi \) to \( \tilde X \) gives a characteristic cover \( \tilde X \to X \), so it follows from Proposition~\ref{prop: rephrase finite type} that \( \tilde X_\circ \) is one ended. In particular, $\tilde X_\circ \smallsetminus K$ has only one component which is not precompact. On the other hand, by Lemma~\ref{lem:noncompact_component}, no component of \( \tilde{X} \smallsetminus \tilde{Y} \) is precompact, so the same is true of $\tilde{X}_\circ \smallsetminus K$. This shows the complement of $K$ in $\tilde{X}_\circ$, and hence in $\tilde{S}_\circ$, is connected. 
\end{proof}

  \subsection{The universal abelian cover}
\label{sec:UAC}

 If $S$ is an orientable borderless surface with fundamental group $\pi_1 (S)$, one example of a characteristic cover of $S$ is the \emph{universal abelian cover} $S^{ab}$, which is defined to be the cover corresponding to the commutator subgroup $[\pi_1(S),\pi_1(S)]$. 
Unless $S$ is a disk or sphere, the universal abelian cover $S^{ab}$ has infinite degree and is not locally finite, so Theorem~\ref{thm: not locally finite} implies that $S^{ab}$ is homeomorphic to either the disk, the flute surface, the Loch Ness monster surface, or the spotted Loch Ness monster surface. 
Given a surface $S$, the next theorem tells you \emph{which} of these surfaces is its universal abelian cover. 
  
  \begin{named}{Proposition \ref{Prop:UAC}}[Universal abelian covers] 
    For a surface $S$, let $S^{ab}$ be its universal abelian cover. Then
    \begin{enumerate}
      \item if $S$ is $\BR^2$, the annulus, or the torus, then $S^{ab}\cong \BR^2$,
      \item if $S$ is the sphere, then so is $S^{ab}$,
      \item if $S$ is the once-punctured torus, then $S^{ab}$ is the flute surface,
      \item if $S$ is a finite-type surface with one puncture and genus at least two, then $S^{ab}$
        is the spotted Loch Ness monster surface,
      \item otherwise, $S^{ab}$ is the Loch Ness monster surface. 
    \end{enumerate}
  \end{named}
  
\begin{proof} Cases (1) and (2) are trivial, so let us assume from now on that $S$ has non-abelian fundamental group.

First, assume that \( S \) is a finite-type genus \( g \) surface with \( n \) punctures. 
We can then realize \( S \) as the interior of a compact surface \( \bar S \) with \( n \) boundary components. 
The abelianization of \( \pi_1 \bar S \) is \( \BZ^{2g+n-1} \), so \( \BZ^{2g+n-1} \) acts properly and cocompactly on \( \bar S^{ab} \); hence, the end space of \( \BZ^{2g+n-1} \) and the end space of \( \bar S^{ab} \) are homeomorphic. As we are assuming that \( \pi_1(S) \) is not abelian, we have that \( 2g+n-1 \geq 2 \), so that \( \BZ^{2g+n-1} \), and hence \( \bar S^{ab} \), is one ended. 

If \( n = 0 \), then \( \bar S = S \), and the classification of surfaces implies that \( S^{ab} \) is the Loch Ness monster surface. 
If $n>1$, then the preimage of $\partial \bar S$ in $\bar S^{ab}$ is a collection of properly embedded lines. Deleting these lines cannot increase the number of ends, so $S^{ab}$ is also one-ended.
Hence, $S^{ab}$ is the Loch Ness monster surface.

Now, suppose $n=1$, which implies that \( g \geq 1 \).
The single boundary component of $\bar S$ lifts homeomorphically, and so $\partial \bar S^{ab}$ is an infinite collection of circles. 
Gluing on a disk to each of these circles results in a borderless one-ended surface \( \tilde S_{cap} \) and a covering \( \tilde S_{cap} \to S_{cap} \) with deck group \( \BZ^{2g} \), where  $S_{cap}$ is the closed surface obtained by gluing on a disk to $\bar S$.
By the classification of surfaces, \( \tilde S_{cap} \) is either the plane or the Loch Ness monster surface.
If it is the plane, then this covering map is the universal covering map and the fundamental group of the surface is free abelian; hence, \( S_{cap} \) is a torus, so \( S \) is the once-punctured torus and its universal abelian cover is obtained by puncturing the plane in an infinite discrete set of points, implying \( S^{ab} \) is the flute surface.
Otherwise, \( \tilde S_{cap} \) is the Loch Ness monster surface, and \( S^{ab} \) is obtained by puncturing the Loch Ness monster surface in an infinite discrete set of points, implying that \( S^{ab} \) is the spotted Loch Ness monster surface.

 Finally,  assume that $S$ is an infinite-type surface. 
 The universal abelian cover $S^{ab}\to S$ is not locally finite, as the deck group is free abelian and hence not a torsion group. 
 By Theorem \ref{thm: not locally finite} and the fact that $\pi_1(S^{ab})$ is nontrivial, the only options are the flute surface, the Loch Ness monster surface, or the spotted Loch Ness monster surface. 
 We claim there are no isolated planar ends in $S^{ab}$. 
 Indeed, regularity of the cover implies that any such end covers an isolated planar end of $S$. 
 But any peripheral curve in $S$ is part of a free basis for $\pi_1(S)$, and hence no power of it lies in the commutator subgroup, so the preimage of any  annular neighborhood of an isolated planar end of $S$ has simply connected components in $S^{ab}$. 
 Hence, $S^{ab}$ is the Loch Ness monster surface.
 \end{proof}

\subsection{$\BZ/n\BZ$-homology covers}
\label{thm: finite homology coverssec}

Recall from the introduction that the \emph{$\BZ/n\BZ$-homology cover} of a surface $S$ is the cover corresponding to the kernel of the map $\pi_1 (S) \to H_1(S,\BZ/n\BZ)$. 

\begin{named}{Theorem \ref{thm: finite homology covers}}[$\BZ/n\BZ$-homology covers]
    Say $S$ is an infinite-type orientable surface and $\pi: \tilde S \to S$ is its $\BZ/n\BZ$-homology cover, where $n\geq 2$. 
    If $S$ has no isolated planar ends, then $\tilde S $ is homeomorphic to the Loch Ness monster surface. Otherwise, $\tilde S$ is homeomorphic to the spotted Loch Ness monster surface.
\end{named}

The main point is to show the following lemma:

\begin{lemma}\label{atmostone!}
With $\pi : \tilde S \to S$ as above, the cover $\tilde S$ has at most one end that is not an isolated planar end.
\end{lemma} 

Deferring the proof of Lemma \ref{atmostone!}, we now prove Theorem \ref{thm: finite homology covers}.
\begin{proof}[Proof Theorem \ref{thm: finite homology covers}]
We first show that 
$\tilde S $ has infinite genus.
Since $\pi$ is a regular cover of infinite degree, it suffices to prove that $\tilde S$ has positive genus. 

If $S$ has positive genus, we can pick an embedded once punctured torus $T \subset S$. Since $S$ has infinite-type, $S \smallsetminus T$ is noncompact, so Mayer-Vietoris implies that $H_1(T,\BZ) \hookrightarrow H_1(S,\BZ)$. Hence, picking a component 
\[\tilde T \subset \pi^{-1}(T) \subset \tilde S~,\]
the cover $\tilde T \to T$ is the $\BZ/n\BZ$-homology cover of $T$. As $\pi_1 (T) \cong F_2$, we have 
\[H_1(T,\BZ/n\BZ) \cong (\BZ/n\BZ)^2~,\] 
so the cover $\tilde T \to T$ has degree $n^2$, so  $\chi(\tilde T) = -n^2$. However, since $\partial T$ is trivial in homology, it lifts homeomorphically to $\tilde T$, and hence $\tilde T$ has $n^2$ boundary components, implying $\tilde T$ is a $n^2$-punctured torus\footnote{One can also just visualize the cover directly by unwrapping the meridian and longitude direction $n$-times.}, so $\tilde S$ has positive genus.

If $S$ has genus zero, take some essential $4$-punctured sphere $T\subset S$. Then $S \smallsetminus T$ has four components. Moreover, all four components of $S \smallsetminus T$ are noncompact: any compact component is a planar surface with one boundary component, and hence is a disk, contradicting that $T$ is essential. It follows that $H_1(T,\BZ) \hookrightarrow H_1(S,\BZ)$. Hence, if $\tilde T \subset \pi^{-1}(T)$ is a component, then $\tilde T \to T$ is the $\BZ/n\BZ$-homology covering of $T$, which has degree 
$n^3$, so $\chi(\tilde T) =-2n^3$. Each component of $\partial T$ is primitive in $H_1(S,\BZ)$, so unwraps $n$ times in $\tilde T$, and hence there are $n^3/n = n^2$ lifts of each, for a total of $4n^2$ boundary components. A quick computation shows that $\tilde T$ has genus $g= 1 - \frac 12(-2n^3+4n^2)$, which is never zero for $n\geq 2$. Thus, by considering the images of $\tilde T$ under the infinite deck group, $\tilde S$ must have infinite genus.

\medskip

We now know that $\tilde S$ has infinite genus. By Lemma \ref{atmostone!}, $\tilde S$ is either the Loch Ness monster surface or its spotted version. If $S$ has isolated planar ends, these ends lift to isolated planar ends in $\tilde S$, so $\tilde S $ is the spotted Loch Ness monster surface. Conversely, since the cover is regular, any isolated planar end in $\tilde S$ covers an isolated planar end in $S$. So if $S$ has no isolated planar ends, $\tilde S$ is the Loch Ness monster surface.
\end{proof}

\subsubsection{The proof of Lemma \ref{atmostone!}}
\label{sec: pflem}
 Before formally proving Lemma~\ref{atmostone!}, we give a brief heuristic of the proof.

The rough idea here is to use the fact that the cover $\tilde S \to S$ is abelian, and the fact that homologically essential curves are ubiquitous in \( S \), and such curves do not lift homeomorphically to $\tilde S$. 
Slightly more precisely, and assuming for simplicity that neither $S$ nor $\tilde S$ has isolated planar ends, we take a compact set $\tilde K \subset \tilde S$ that projects to $K\subset S$, a path $\tilde a $ in $\tilde S$ with endpoints outside of $\tilde K$, and we claim that there  exists another path $\tilde \gamma$ with the same endpoints that does not enter $\tilde K$ (this will show that $\tilde S$ is one-ended). 
To construct such a $\tilde\gamma$ we work as follows. 
Writing $a = \pi(\tilde a)$ and $\gamma = \pi(\tilde \gamma)$, it may be that there is not really another way to travel within $S$ between the endpoints of $a$ without roughly following $a$, which might take us through $K$. 
So, we try to construct $\tilde \gamma$ so that at some point its projection $\gamma$ follows $a$, but when it does so, $\tilde \gamma$ still avoids going through $\tilde K$. 
To do this, we first make $\gamma$ loop around some homologically essential curve $b$, so that when $\gamma$ tracks $a$ the lift $\tilde \gamma$ is actually tracking a translate of $\tilde a$ rather than $\tilde a$ itself. 
We then loop backwards around $b$ at the end, and use that the cover is abelian to say that this cancels out the first loop around $b$, so that $\tilde \gamma$ terminates at an endpoint of $\tilde a$ rather than at some translate of that endpoint.

\medskip

We now prove Lemma~\ref{atmostone!} by contradiction. 
Suppose that $\tilde S$ has more than one end that is not isolated and planar. 
We first show:

\begin{claim}
    There exists a connected, finite-type subsurface $X\subset S$ such that 
    \begin{enumerate}
        \item every component of $S\smallsetminus X$ has infinite type, and
        \item  there is a finite-type connected component $\tilde X$ of $\pi^{-1}(X) \subset \tilde S$ such that the complement $\tilde S \smallsetminus \tilde X$ contains two distinct infinite-type connected components $\tilde U_+$ and $\tilde U_-$. 
    \end{enumerate}
\end{claim}
\begin{proof}
    Since $\tilde S$ has more than one end that is not isolated and planar, there is a compact subsurface $\tilde K \subset \tilde S$ that has two distinct infinite-type connected components. 
The projection $\pi(\tilde K) \subset S$ is compact, so there is a finite-type subsurface $X \subset S$ containing the projection. 
We can assume that all components of $S \smallsetminus X$ have infinite type, by adding any finite type components into $X$. 
Since $\pi$ is locally finite (see \S \ref{locallyfinite}), all components of $ \pi^{-1}(X)$ have finite type, and the component $\tilde X$ containing $\tilde K$ has at least two infinite-type complementary components.
\end{proof}

Pick an oriented path $\tilde a$ in $\tilde X$ that starts at a point $\tilde x_- \in \partial \tilde U_-$ and ends at a point $\tilde x_+\in \partial \tilde U_+$ and let $a,x_-,x_+$ be the projections to $S$. 
Then $x_\pm$ lies in $\partial X$ and is adjacent to a component $V_\pm \subset S \smallsetminus X$. 
As $V_\pm$ has infinite type, there is a nonseparating bi-infinite arc $\alpha_\pm$ in $V_\pm$. 
If $V_+=V_-$, we can choose the arcs $\alpha_+,\alpha_-$ so that the union $\alpha_+\cup \alpha_-$ is nonseparating in $V_+=V_-$. 
Let $ b_\pm$ be an oriented loop in $V_\pm$ that starts and ends at $x_\pm$, and that intersects $\alpha_\pm$ exactly once. 
If $V_+=V_-$, we can assume that $b_\pm$ is disjoint from $\alpha_\mp$.

\begin{figure}[t]
\centering
	\def\svgwidth{0.85\textwidth}
 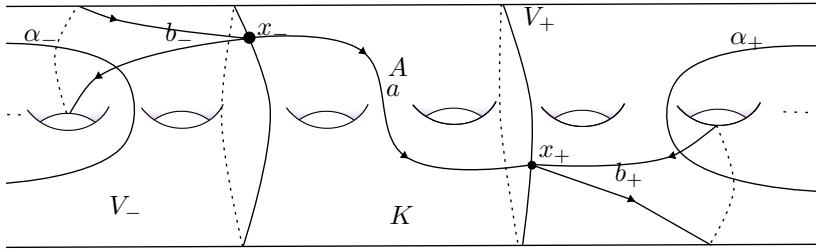
	\caption{The compact subsurface $K$, the reference bi-infinite arcs \( \alpha_- \) and \( \alpha_+ \), and the arcs \( b_- \), \( b_+ \), and \( a \) used to form \( \gamma \).}
	\label{fig:gammaloop}
\end{figure}

Consider now the path $\gamma = b_- a b_+^{-1} a^{-1} b_-^{-1} a b_+$ in $S$, which starts at $x_-$ and ends at $x_+$, see Figure \ref{fig:gammaloop} for a configuration.

Then $\gamma a^{-1}$ is a loop that starts and ends at $x_-$, and since $$\gamma a^{-1}= b_- (a b_+^{-1} a^{-1}) b_-^{-1} (a b_+ a^{-1}) = [b_-,ab_+^{-1}a^{-1}],$$
this loop is trivial in homology and therefore lifts to $\tilde S$. 
Choose the lift such that the terminal $a^{-1}$ lifts to the already defined path $\tilde a^{-1}$, which starts at $\tilde x_+$ and ends at $\tilde x_-$. 
If we do this, then the resulting lift $\tilde \gamma$ of $\gamma$ starts at $\tilde x_-$ and ends at $\tilde x_+$.

\begin{claim}\label{disjointclaim}
    $\tilde \gamma$ is disjoint from $\int(\tilde X)$.
\end{claim}

This claim will finish the proof of Lemma~\ref{atmostone!}, as it contradicts the fact that $\tilde x_-$ and $\tilde x_+$ lie on the boundaries of distinct connected components of $\tilde S \smallsetminus \tilde X$. 
So, all that is left is to prove the claim.

\begin{proof}[Proof of Claim  \ref{disjointclaim}]
Below, when we talk about `the lift of' a path like $b_-, $ or $a,$ or $ b_+^{-1}$, etc.,~we mean the lift that appears as the associated subpath of $\tilde \gamma$. 
Since there are two traversals of $a$ in $\gamma$, we refer to the two lifts as the `first lift' of $a$ and the `second lift' of $a$, using the order they appear in the word. 
In this terminology, orientations matter: the `lift of $b_-^{-1}$' may not be obtained by traversing the `lift of $b_-$' backwards. 

First, observe that all the lifts of $b_-,b_-^{-1},b_+,b_+^{-1}$ are disjoint from $\int(\tilde X)$, as their projections are disjoint from $X$ by construction.
Now, suppose the first lift of $a$ intersects $\int(\tilde X)$. 
Then we can create a loop in $\tilde S$ by concatenating the lift of $b_-$, part of the lift of $a$, and a path in $\tilde X$. 
The projection of this loop to $S$ intersects $\alpha_-$ exactly once, as the latter two parts project into $X$. 
This is impossible, since the projection must be trivial in $\BZ/n\BZ$-homology.
The same argument shows that the lift of $a^{-1}$ is disjoint from $\int(\tilde X)$: if not, one can construct a loop that is trivial in $\BZ/n\BZ$-homology and intersects $\alpha_-$ once.

A similar argument shows that the second lift of $a$ is disjoint from $\int(\tilde X)$: if not, one can construct a loop that is trivial in $\BZ/n\BZ$-homology and intersects $\alpha_+$ once. 
Note that above, we could have used $\alpha_-$ or $\alpha_+$, but here we need to use $\alpha_+$, as both $b_-$ and $b_-^{-1}$ appear before the second occurrence of $a$ in the word defining $\gamma$.
\end{proof}

\section{Arbitrary covers}\label{sec: arbitrary covers}
In this section we prove the following theorem.
\begin{named}{Proposition \ref{Prop:arbitrary covers}}[Everything covers everything]
 Suppose that $S$ is an orientable, borderless surface with non-abelian fundamental group. Then $S$ is covered by any noncompact borderless orientable surface. 
\end{named} 
The idea for the proof of Proposition \ref{Prop:arbitrary covers} is to show that any surface $S$ with non-abelian fundamental group is covered by the blooming Cantor tree surface $\BCT$ and that any infinite-type surface admits a $\pi_1$-injective embedding into $\BCT$, see Proposition \ref{prop:embedsurf}. Then, Lemma \ref{lem:coverofsubsurf} gives us Proposition \ref{Prop:arbitrary covers}.

Before starting the proof of Proposition \ref{Prop:arbitrary covers} we show that the blooming Cantor tree surface regularly covers the Loch Ness monster surface.

\begin{lemma}\label{bctlnm} The blooming Cantor tree surface regularly covers the Loch Ness monster surface.
\end{lemma} 
\begin{proof} Let $G$ be the Cayley graph of $\mathbb Z^2$. Note that if we replace every vertex of $G$ by a four-punctured torus and we glue them according to the edges relation we obtain a surface $L$ that, by classification, is homeomorphic to the Loch Ness monster surface. Let $\tilde G$ be the universal cover of $G$, i.e., $\tilde G$ is a complete 4-valent tree. Then, if we replace each vertex of $\tilde G$ by a four punctured torus we obtain a surface $\tilde S$ that, again by classification, is homeomorphic to the blooming Cantor tree surface. Then the induced covering map from the graphs induces the required covering map from $\tilde S$ to $S$.\end{proof}

We now show that any infinite-type surface can be embedded in the blooming Cantor tree surface $\CB$ as a connected component of the complement of a collection of simple closed curves and bi-infinite arcs.

\begin{proposition}\label{prop:embedsurf} For any borderless non-compact orientable surface $S$, there is a $\pi_1$-injective embedding of $S$ into the blooming Cantor tree $\BCT$.
\end{proposition}
\begin{proof} 

We assume throughout the proof that $S$ has either zero or infinite genus, since an arbitrary  surface can be written as the connected sum $S \# T$ of such an $S$ and a closed orientable surface $T$, and a $\pi_1$-injective embedding $S \hookrightarrow \BCT$ induces a $\pi_1$-injective embedding of $S \# T \hookrightarrow \BCT \# T \cong \BCT$. 

Let $(\mathcal E, \mathcal E_g)$ be the genus-marked end space of $S$, which is nonempty since $S$ is non-compact. Since the end space $\CE(\BCT)$ is a Cantor set, we can fix an embedding $\iota: \mathcal E\rightarrow \CE(\BCT)$, and we let $E=\iota(\mathcal E)$ and $E_g=\iota(\mathcal E_g)$. We will construct an open, connected $\pi_1$-injective subsurface $S'' \subset \BCT$ such that
\begin{enumerate}
    \item[(a)] the inclusion $S''\hookrightarrow \BCT$ induces a homeomorphism $(\CE(S''), \CE_g(S'')) \cong (E, E_g)$,
    \item[(b)] the genus of $S''$ is the same as the genus of $S$.
\end{enumerate}
By classification of surfaces, see Theorem~\ref{thm: classification}, such a surface $S''$ will be homeomorphic to $S$, proving the proposition.

As a first step, we construct a $\pi_1$-injective \emph{closed} subsurface $S'\subset \BCT$ (with boundary) that has properties (a) and (b).
Let $\Lambda \subset \BCT$ be a union of simple closed curves that cuts $\BCT$ into a collection of thrice-punctured tori, connected together in the pattern of a trivalent tree. For each of the thrice-punctured tori $C \subset \BCT \setminus \Lambda$, cut $C$ along some simple closed curve into a punctured torus $T_C$, and a four-punctured sphere $X_C \subset C$. See Figure \ref{fig:B} for such a configuration.

\begin{figure}[htb!]
\centering
\includegraphics[width=.8\textwidth]{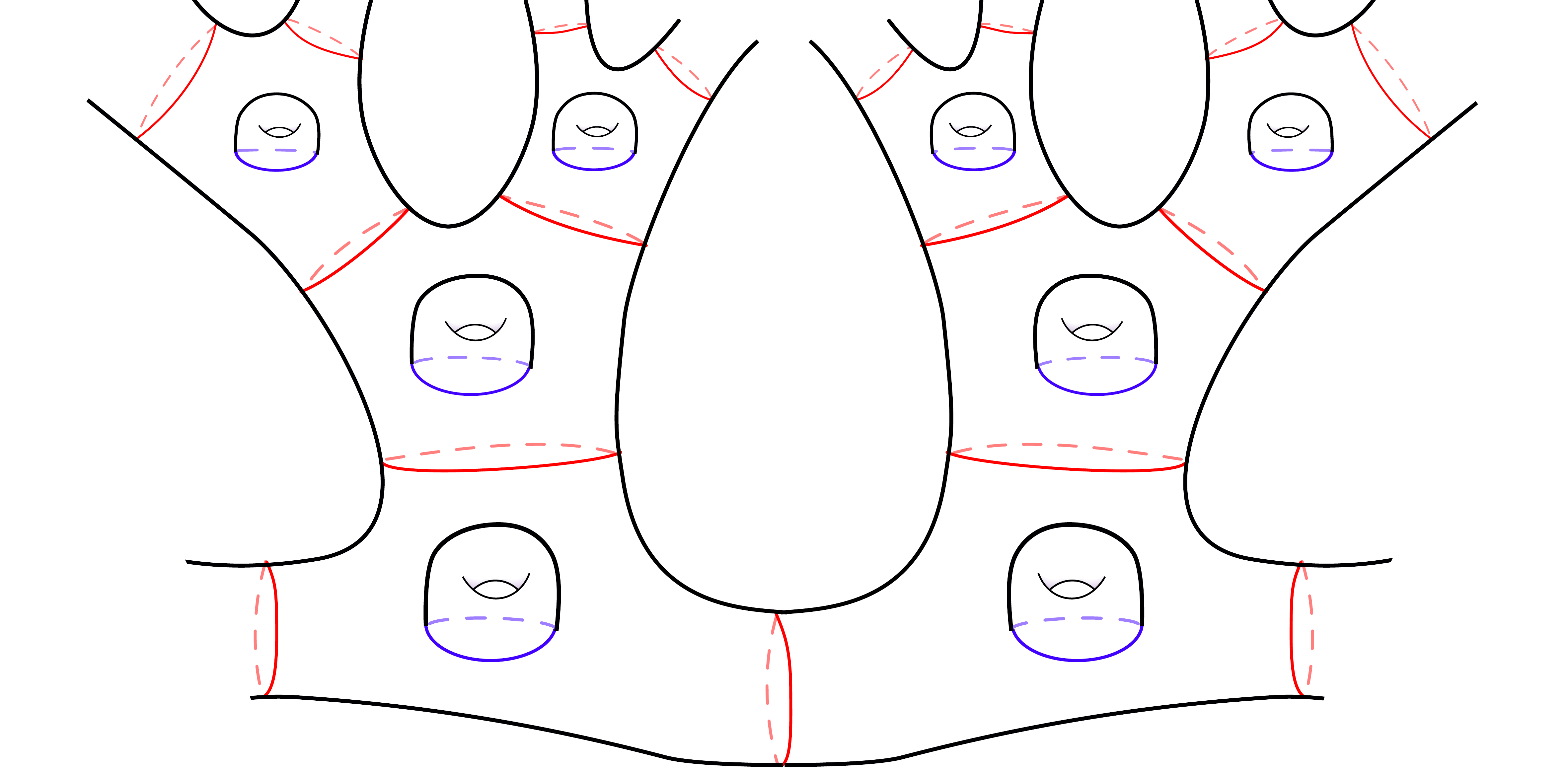}
\caption{The blooming Cantor tree surface with multicurve $\Lambda$ depicted in red.}
\label{fig:B}
\end{figure}

Fixing some  component $\lambda_0 \subset \Lambda$ as the `root', we then define $S' \subset \BCT$ to be the union of:
\begin{enumerate}
    \item all $\overline{X_C}$ such that $C$ separates $\lambda_0$ from some point of $E$,
    \item all $\overline{T_C}$ such that $C$ separates $\lambda_0$ from some point of $E_g$,
\end{enumerate} 
here $\overline{X_C}$ and $\overline{T_C}$ are just the closures of the open subsurfaces mentioned above. Observe that $S'$ is connected, it $\pi_1$-injects into $\CB$, there is a homeomorphism $(\CE(S'), \CE_g(S')) \simeq (E, E_g)$ as above and the genus of $S'$ is the same as that of $S$. The boundary $\partial S'$ consists of countably many circles. Our goal now is to remove the compact boundary components of $S'$ so that they are replaced by half-planes.

We claim that for \emph{any} orientable surface $S'$ such that $\partial S'$ is a union of simple closed curves, there is a $\pi_1$-injective borderless subsurface $S'' \subset S'$, with the same genus, such that the inclusion induces a homeomorphism of genus-marked end spaces. Let $\{\gamma_i\}_{i\in\mathbb N}$ be the collection of boundary components of $S'$ and let $\{r_i\}_{i\in\mathbb N}$ be a locally finite collection of pairwise disjoint properly embedded rays on $S'$ such that $r_i$ starts on $\gamma_i$, and let $P_i$ be a closed regular neighborhood of $\gamma_i \cup r_i$ in $S'$. If we take these regular neighborhoods small enough, they are all disjoint and we obtain a surface $S'' := S' \setminus \cup_i P_i$. Observe that the end space of $S''$ agrees with that of $S'$, so it also satisfies the required properties.

Composing the inclusions $S'' \hookrightarrow S' \hookrightarrow \BCT$ gives an open, $\pi_1$-injective subsurface of $\BCT$ satisfying (a) and (b) above, which proves the proposition. \end{proof}

\begin{figure}[htb!]
\centering

\includegraphics[width=.85\textwidth]{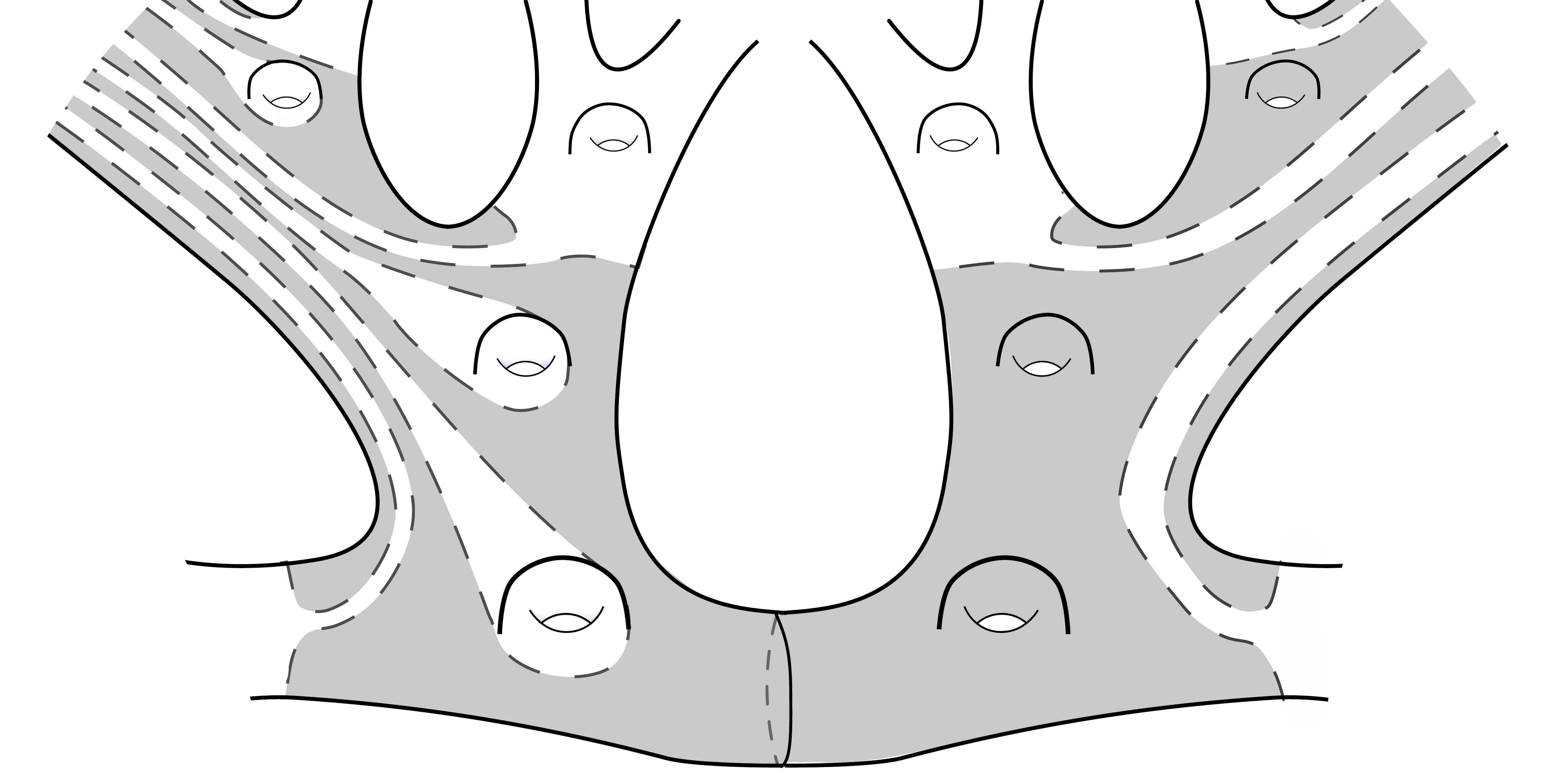}
\caption{The once-punctured Loch Ness monster surface as a subsurface of $\BCT$ according to our construction. The left side of the figure represents the sole planar end, while the right side represents the 
 end accumulated by genus.\label{fig:ray}}
\end{figure}

\begin{proof}[Proof of \Cref{Prop:arbitrary covers}] By Proposition \ref{prop:embedsurf} and Lemma \ref{lem:coverofsubsurf} it suffices to show that given $S$ there exists a covering space $B$ homeomorphic to a blooming Cantor tree surface. Since $S$ has non-abelian $\pi_1$ we can find a $\pi_1$-injective pair of pants $P\subset S$. Then, by Lemma \ref{lem:coverofsubsurf} we get a cover $S''\to S$ where $S''$ is an open pair of pants. By Proposition \ref{Prop:UAC} $S''$ is covered by a Loch Ness monster surface $L$ and we conclude by Lemma \ref{bctlnm}.
\end{proof}

\appendix
\section{Ends, group cohomology, and Swarup's theorem}\label{sec: appendix}
The goal of this appendix is to show how Theorem \ref{Thm:relativestallings} is a restatement of the main result in Swarup's paper \cite{swarup1977relative}. 
To do so, we relate three ways of looking at ends of a group $G$. 
In the background, we have the space of ends $\CE(G)$ as defined in \S \ref{sec: ends}. 
We want to relate this to $H^1(G,\BZ G)$, which is the definition of the space of ends in Swarup's paper \cite{swarup1977relative}. 
To do so, we relate both to a combinatorial description of ends that is similar to a description in Stallings's paper \cite{stallings1971group}. 
All of this material is likely well known to experts, but we did not readily find a reference for the portions relevant  to the goals of this article, and so we present it here in hopes of providing a useful reference to others. 

\subsection{Group cohomology}
Let $G$ be a group. 
A \emph{$G$-module} is an abelian group that comes equipped with an action of $G$ by automorphisms. 
For instance, any abelian group $A$ with the trivial $G$ action is a $G$-module. 
On the other end of the spectrum, let
$$\BZ^G := \{ \phi\co G \to \BZ \}$$
denote the set of functions from $G$ to $\BZ$, which we consider as an abelian group under addition.
Given $\phi \in \BZ^G$ and $g\in G$, let $g\phi$ be the function defined by
\begin{align*} g\phi:\, G& \to \BZ \\
h&\mapsto \phi(g^{-1}h). 
\end{align*}
We realize \( \BZ^G \) as a \( G \)-module via this (left) \( G \)-action. 

The \emph{group ring of \( G \)}, denoted by \( \BZ G \), is the submodule of $\BZ^G$ given by 
$$\BZ G := \{ f\co G \to \BZ \ | \ f \text{ has finite support}\},$$
where the \emph{support} of $f$ is the set of all $g\in G$ with $f(g)\neq 0$. 
We can express the elements of $\BZ^G$  as formal $\BZ$-linear combinations
$$\sum_{g\in G} c_g \cdot g,$$
where \( c_g \in \BZ \) and---abusing notation---we view \( g \) as the indicator function that assigns \( 1 \) to \( g \) and \( 0 \) to all other group elements.
In this form, the left $G$ action on $\BZ^G$ described above is given by 
$$h\cdot\left(\sum_{g\in G} c_g \cdot g\right) = \sum_{g\in G} (c_g \cdot hg)$$
for $h \in G$. 
The elements of $\BZ G$ can then be realized as the formal sums in which $c_g = 0$ for all but finitely many g. 
We will switch back and forth between formal $\BZ$-linear combinations and functions as convenient.

Suppose that $M$ is a $G$-module. 
We will now define the cohomology groups $H^n(G,M)$ with coefficients in $M$. 
An \emph{$n$-cochain with coefficients in \( M \)} is  a function $\phi : G^n \to M$, where $G^n = G \times \cdots \times G$, and by convention $G^0 = \{1\}$. 
Let \( C^n \) denote the set of \( n \)-cochains. 
We then define the coboundary maps $\partial^n : C^n \to C^{n+1}$ as follows:
for \( n > 0 \), set
\begin{align*}
(\partial^n\phi)(g_1,\ldots,g_{n+1}) \ \ \ = \ \ \ \ & g_1 \phi(g_2,\ldots,g_{n+1}) 
\\ + \ & \sum_{i=1}^n(-1)^i \phi(g_1,\ldots,g_ig_{i+1},\ldots,g_{n+1}) 
\\ + \ & (-1)^{n+1} \phi(g_1,\ldots,g_{n}).
\end{align*}
and for \( n = 0 \), set
\begin{equation}
(\partial^0\phi)(g) = g \cdot \phi(1) - \phi(1) = (g-1) \phi(1). \label{partial0}
\end{equation}
The cohomology groups are now defined to be 
$$H^n(G,M) := \mathrm{ker} \ \partial^n / \mathrm{Im} \ \partial^{n-1}.$$
This definition relies on the fact that \( \mathrm{Im}(\partial^{n-1}) \subset \ker\partial^n \), a fact we have not checked; we will verify this for \( n =1 \) below, which is the only case we will use in this appendix. 

Now, $0$-chains are functions $\phi : \{1\} \to M$, or in other words, a choice of an element of \( M \). 
By \eqref{partial0},  $\partial^0\phi =0$ exactly when the element $\phi(1) \in M$ is invariant under the action of $G$, and hence we get
$$H^0(G,M) \cong \{ x\in M \ | \ gx=x, \forall g\in G\}.$$
Here are some examples of different $G$-modules $M$.
\begin{enumerate}
\item If $M$ is an abelian group with the trivial $G$-action, then $H^0(G,M) \cong M$.
\item If $M=\BZ^G$, then $H^0(G,M)\cong \BZ$ since the $G$-invariant elements of $\BZ^G$ are exactly those where all the coefficients are the same.
\item If $M=\BZ G$, then $$H^0(G,\BZ G) \cong \begin{cases}
 \BZ & \text{ if } G \text{ is finite } \\ 
0 & \text{ if } G \text{ is infinite. }
 \end{cases}$$
 Indeed, if \( [f] \in H^0(G, \BZ G) = \ker \partial^0 \), then the support of \( f \) is finite and \( f(g) = f(h) \) for all \( g, h \in G \). 
\end{enumerate}

To compute $H^1(G,M)$, we need to understand the image of $\partial^0$ (i.e., the set of \emph{coboundaries}) and the kernel of $\partial^1$ (i.e., the set of \emph{cocycles}). 
Given \( x \in M \), define the 1-cochain \( \phi_x \co G \to M \) by 
$$ \phi_x(g) = (g-1)x.$$
In \eqref{partial0}, we computed the image of \( \partial^0 \), establishing that the set of 1-coboundaries is \( \{ \phi_x : x \in M \} \). 

To understand $1$-cocycles, first note that if $\phi $ is a $1$-cochain, then
$$(\partial^1\phi)(g,h) = g\phi(h)-\phi(gh)+\phi(g).$$
Setting the right-hand side to zero, we see that $\phi$ is a $1$-cocycle if and only if 
$$\phi(gh) = g \phi(h) + \phi(g), $$
i.e., \( \phi \) satisfies the \emph{cocycle condition}. 
Let us verify that every 1-coboundary is a 1-cocycle, i.e., \( \mathrm{Im}(\partial^0) \subset \ker \partial^1 \). 
Fix a 1-cobounday \( \phi_x \), where \( x\in M \). 
Then 
\begin{equation}
\phi_x(gh) = (gh-1)x = g(h-1)x + (g-1)x = g\phi_x(h)+\phi_x(g) \label{sq}
\end{equation}
implying that \( \phi_x \) satisfies the cocycle condition and is therefore a 1-cocycle. 

Let us work out two example computations of \( H^1(G, M) \). 
\begin{enumerate}
\item[(A)] Suppose that $M$ is an abelian group with a trivial $G$-action. 
In this case, the cocycle condition becomes $\phi(gh) = \phi(g) + \phi(h)$, i.e.\ $\phi$ is a homomorphism. 
Given a coboundary \( \phi \), there exists \( x \in M \) such that \( \phi = \phi_x \). 
It follows that \( \phi(g) = (g-1)x = gx - x = x -x = 0 \), and hence \( \phi = 0 \).
As all coboundaries vanish, 
$$H^1(G,M) = \mathrm{Hom}(G,M).$$
In the case \( M = \BZ \), we have that $H^1(G,\BZ)$ is isomorphic to the free part of the abelianization of $G$.

\item[(B)] Suppose that $M=\BZ^G$. We claim that $H^1(G,\BZ^G)=0$. To see this, suppose \( \phi\co G \to \BZ^G \) is a 1-cocycle.
Define \( x \in \BZ^G \) by 
\[ x(g) = \phi(g^{-1})(1). \]
We claim that \( \phi = \phi_x \), and hence \( \phi \) is a coboundary. 
Indeed, if \( h \in G \), then:
\begin{align*}
\phi_x(g)(h)	&= (g-1)x(h) \\
			&= gx(h) - x(h) \\
			&= x(g^{-1}h)-x(h) \\
			&= \phi(h^{-1}g)(1)-\phi(h^{-1})(1) \\
			&= (\phi(h^{-1}g)-\phi(h^{-1}))(1) \\
			&= h^{-1}\phi(g)(1) \\
			&= \phi(g)(h)
\end{align*}
where the third equality uses the definition of the action and the second-to-last equality uses the cocycle condition.
Therefore, \( \phi = \phi_x \), implying that \( \phi \) is a coboundary.
It follows that \( H^1(G, \BZ^G) = 0 \). 

\end{enumerate}

\subsection{Ends via group cohomology}

In this section we describe how the set of ends of $G$ is encoded in the cohomology group $H^1(G,\BZ G)$. 
It is well known that the number of ends of an infinite finitely generated group \( G \) is equal to \( 1+\mathrm{rank}(H^1(G, \BZ G)) \), where rank denotes the cardinality of a maximal \( \BZ \)-linearly independent subset (see for example \cite[Theorem~13.5.5]{geoghegan2007topological}). 
However, we require a more transparent proof of this fact that will allow us to more easily generalize to the relative setting discussed in Theorem~\ref{Thm:relativestallings}.

\medskip

Suppose $G$ is generated by a finite set $S$ that is symmetric under inversion. Let ${Cay}(G)$ be the associated \emph{left Cayley graph} for $G$, that is, the graph whose vertex set is \( G \) and where \( g \) and \( h \) are adjacent if \( h = sg \) for some \( s \in S \). 
Note that this is opposite the usual construction of Cayley graphs; however, we need to use left Cayley graphs here since we chose to consider $\BZ G$ as a left $G$-module. 

Given $x\in \BZ^G$, the \emph{boundary} of $x$ is the subset $\partial x \subset G$ consisting of all $g\in G$ such that $x(sg)\neq x(g)$ for some $s \in S$.  
For example, if $G=\BZ$, \( S = \{\pm 1\} \), and
\[
x(n) = \left\{
\begin{array}{ll}
1 &  \text{if } n \geq 0 \\
0 &  \text{if } n < 0 
\end{array}\right.
\]
then $\partial x = \{-1,0\}.$

Recall that, given \( x \in \BZ^G \), we defined \( \phi_x \co G \to \BZ^G \) by \[ \phi_x(g) = (g-1)x .\]
The next two lemmas tell us that if \( x \) has finite boundary, then \( \phi_x \) is a 1-cocycle with \( \BZ G \) coefficients and that these account for all 1-cocycles with \( \BZ G \) coefficients. 
The third lemma characterizes the 1-coboundaries with \( \BZ G \) coefficients.

\begin{lemma}\label{finiteimage}
If \( x \in \BZ^G \), then \( \phi_x \) is a 1-cochain with \( \BZ G \) coefficients if and only if \( \partial x \) is finite. 
\end{lemma}

\begin{proof}
If \( g \in \partial x \), then there exists \( s \in S \) such that \( x(sg) \neq x(g) \).
It follows that \( s^{-1}x(g) \neq x(g) \) and hence \( s^{-1}x(g) - x(g) \neq 0 \), which says that \( g \) is an element in the support of \( (s^{-1}-1)x \). 
It follows that \( \partial x \) is contained in the union of the supports of the elements \( (s-1)x \in \BZ^G \) for \( s \in S \). 
If we assume that \( \phi_x \) is a 1-cochain with \( \BZ G \) coefficients, then \( (s-1)x = \phi_x(s) \in \BZ G \) for all \( s \in S \); it follows that \( \partial x \) is finite, establishing the forwards direction.  

Now, assume $\partial x$ is finite. 
Fix \( g \in G \).
If $g^{-1}=s_1\cdots s_m$, where $s_i \in S$, we claim that $(g-1)x$ is supported within the $m$-neighborhood of $\partial x \subset G$, with respect to the word metric. Indeed, if $h$ lies outside this neighborhood, then for each $i$, the elements $s_i\cdots s_m h$ does not lie in $\partial x$, so $x(s_i\cdots s_m h) = x(s_{i-1}\cdots s_m h)$, and inductively we get $x(g^{-1}h)=x(h)$.
As \( x(g^{-1}h) = gx(h) \), we have \( gx(h)-x(h) = 0 \), so $(g-1)x(h)=0$ implying that \( \phi_x(g) = (g-1)x \in \BZ G \).
As \( g \in G \) was arbitrary, we conclude that \( \phi_x \) is a 1-cochain with \( \BZ G \) coefficients. 
\end{proof}

\begin{lemma}[$1$-cocycles with $\BZ G$-coefficients] 
A 1-cochain \( \phi \) with \( \BZ G \) coefficients is a 1-cocycle if and only if there exists \( x \in \BZ^G \) with finite boundary such that \( \phi = \phi_x \). 
\label{cocycles}
\end{lemma}

\begin{proof}
Let \( x \in \BZ^G \) with finite boundary. 
Lemma~\ref{finiteimage} tells us that \( \phi_x \) is a 1-cochain with \( \BZ G \) coefficients, and then \eqref{partial0} tells us that \( \phi_x \) is a 1-cocycle with \( \BZ G \) coefficients. 
Conversely, if \( \phi \) is a 1-cocycle with \( \BZ G \) coefficients, then the computation in (B) above tells us that \( \phi = \phi_x \) with \( x \in \BZ^G \) defined by \( x(g) = \phi(g^{-1})(1) \). 
It is left to check that \( x \) has finite boundary:
as \( \phi_x = \phi \), we have that \( (g-1)x = \phi(g) \in \BZ G \) for all \( g \in G \); hence, by Lemma~\ref{finiteimage}, \( x \) has finite boundary. 
\end{proof}

\begin{lemma}
\label{lem: constant}
Let \( x \in \BZ^G \).
If \( \phi_x \) is a 1-cocycle with \( \BZ G \) coefficients, then \( \phi_x \) is a 1-coboundary with \( \BZ G \) coefficients if and only if \( x \) is constant outside a finite set of \( G \). 
\end{lemma}

\begin{proof}
By definition, \( \phi_x \) is a 1-coboundary if and only if there exists \( y \in \BZ G \) such that \( \phi_x = \phi_y \). 
For the forwards direction, suppose that \( \phi_x \) is a 1-coboundary.
Then there exists \( y \in \BZ G \) such that \( (g-1)x = (g-1)y \) for every \( g \in G \). 
In particular, \( g(x-y) = x-y \) for each \( g \in G \), implying that \( x-y \) is constant.
Hence, $x$  is constant off the finite support of $y$.
Conversely, suppose that \( x \) is constant outside a finite subset \( Y \) of \( G \). 
Let \( d \in \BZ \) such that \( x(g) = d \) for all \( g \notin Y \). 
Define \( y \in \BZ G \) by \( y(g) = 0 \) if \( g \notin Y \) and \( y(g) = x(g) - d \) if \( g \in Y \). 
Then \( x - y \) is constant, and by repeating the computation above in reverse, \( \phi_x = \phi_y \), implying that \( \phi_x \) is a 1-coboundary 
\end{proof}

Now, suppose that $G$ is an infinite finitely generated group.
Let $\CE(G)$ be the space of ends of $Cay(G)$, and equip $Cay(G)\cup\CE(G)$ with the usual topology. 
Identifying \( G \) with the vertex set of \( Cay(G) \), we can view \( G \cup \CE(G) \) as a topological space. 
It is an exercise to show that \( \CE(G) \) is either a singleton, doubleton, or a perfect set; we define \( e(G) \) to be one, two, or \( \infty \), respectively. 

If $x\in \BZ^G$ has finite boundary, then any end of $G$ has a neighborhood in $G$ on which $x$ is constant, so there is a unique function $ \bar x \co \CE(G) \to \BZ$ such that $x\cup \bar x \co G \cup \CE(G) \to \BZ$ is continuous.

Let $C(\CE(G),\BZ)$ be the set of all continuous functions  $\CE(G)\to \BZ$, which is an abelian group under addition. 
Let $Z \subset C(\CE(G),\BZ)$ be the constant functions.
As mentioned above, the rank of an abelian group \( A \), denoted \( \mathrm{rank}(A) \), is the cardinality of a maximal \( \BZ \)-linearly independent subset of \( A \).

\begin{theorem}\label{bigthm}
Let \( G \) be an infinite finitely generated group.
 The map \[ H^1(G,\BZ G) \to C(\CE(G),\BZ)/Z \] defined by $[\phi] \longmapsto [\bar x]$ is an isomorphism, where \( x \in \BZ^G \) satisfies \( \phi = \phi_x \). Consequently, we have 
\[ e(G) = 1+\mathrm{rank}(H^1(G, \BZ G)). \] 
\end{theorem}

Note that we are supposing $G$ is infinite above, so that $\CE(G)$ is nonempty. 
If $G$ is finite, then \( e(G) = 0 \) and, by Lemma~\ref{cocycles}, $H^1(G,\BZ G) = 0$.

\begin{proof}[Proof of Theorem~\ref{bigthm}]
First, we claim the map is well defined: suppose \( [\phi] = [\psi] \).
Then there exists \( x,y \in \BZ^G \) such that \( \phi = \phi_x \), \( \psi = \phi_y \), and, by Lemma~\ref{lem: constant}, \( x-y \) is constant, implying \( [\bar x] = [\bar y] \), as desired. 
The map is a homomorphism as 
$$[\phi_x]+[\phi_y] = [\phi_x+\phi_y] = [\phi_{x+y}]$$
and the map $x \to \bar x$ is linear. 

Next, we show the map is injective: suppose \( x \in \BZ^G \) with finite boundary such that \( \bar x \in Z \).
As \( \bar x \in Z \), there exists \( d \in \BZ \) such that \( \bar x(\xi) = d \) for all \( \xi \in \CE(G) \). 
It follows that \( x(g) = d \) for every \( g \) that lies in an infinite component of \( Cay(G) \smallsetminus \partial x \).
As \( \partial x \) is finite, Lemma~\ref{lem: constant} implies that \( [\phi_x] = 0 \). 
Now, we turn to surjectivity. 
If \( f \in C(\CE(G), \BZ) \), then the continuity of \( f \) and the compactness of \( \CE(G) \) imply that the image of \( f \) is finite; let \( \{d_1,\ldots,d_m\} \subset \BZ \) be the image of \( f \). 
The preimages $E_i:=f^{-1}(d_i)$ are clopen sets that partition $\CE(G)$; therefore, there exists a finite subset \( F \) of \( G \) such that each of the \( E_i \) are the accumulation points in \( \CE(G) \) of a component \( C_i \) of \( Cay(G) \smallsetminus F \). 
Define $x\in \BZ^G$ such that $x(g)=d_i$ for all $g\in C_i$ and $x(g)=0$ otherwise. 
The boundary $\partial x$ is contained in the $1$-neighborhood of $F$, so \( \partial x \) is finite, and $f=\bar x$.
By Lemma~\ref{cocycles}, \( \phi_x \) is a 1-cocycle, and by construction, the image of \( [\phi_x] \) is \( [f] \), establishing the surjectivity of the map. 

We have established that the given map is an isomorphism, so it is left to check the that \( e(G) = 1+\mathrm{rank}(H^1(G, \BZ G)) \). 
First suppose that \( |\CE(G)| = n < \infty \).
Let \( \CE(G) = \{ \xi_1, \ldots, \xi_n \} \) and define \( f_i \in C(\CE(G), \BZ) \) by \( f_i(\xi_i) = 1 \) and \( f_i(\xi_j) = 0 \) if \( j \neq i \). 
Then \( C(\CE(G), \BZ) \) is generated by \( \{ f_1, \ldots, f_n\} \), and it is clear that this set of generators is \( \BZ \)-linearly independent; in particular, \( C(\CE(G), \BZ) \cong \BZ^n \).
Therefore, the quotient \( C(\CE(G), \BZ)/ Z \cong \BZ^{n-1} \), yielding the desired result. 
Now, suppose that \( \CE(G) \) is infinite. 
For any $n \in \BN$, there exists pairwise-disjoint clopen sets \( P_1, \ldots, P_{n+1} \) of \( \CE(G) \) such that $\CE(G) = P_1 \cup P_2 \cup \cdots \cup P_{n+1}$. 
Arguing as in the finite end case, the functions \( f_i \in C(\CE(G), \BZ) \) defined to be 1 on \( P_i \) and 0 on the complement of \( P_i \) generate a copy of \( \BZ^{n+1} \) in \( C(\CE(G), \BZ) \), and therefore their equivalence classes generate a copy of \( \BZ^n \) in \( C(\CE(G), \BZ)/ Z \). 
It follows that the rank of \( C(\CE(G), \BZ) \) is infinite, yielding the desired result. 
\end{proof}

It not relevant to the work that follows, but we note that as \( H^1(G,\BZ G) \) is a free abelian group (see \cite[Theorem~13.5.3]{geoghegan2007topological}), if \( \CE(G) \) is infinite, then \( H^1(G, \BZ G) \) is isomorphic to a countably infinite direct sum of copies of \( \BZ \).

We can now state Stallings's theorem in terms of group cohomology.

\begin{theorem}
If $G$ is a finitely generated group such that $H^1(G,\BZ G)$ is nontrivial, then $G$ admits an action on a simplicial tree $T$ with finite edge stabilizers, no edge inversions and no global fixed point.
\qed
\end{theorem}

  \subsection{Swarup's Theorem}
  Swarup \cite{swarup1977relative} has proved a relative version of Stallings's Theorem, which we now describe. Let $G$ be a finitely generated group, and let $H\subset G$ be a subgroup. 
There is a natural \emph{restriction map} \( r : H^1(G,\BZ G) \to H^1(H, \BZ G) \) given by  
\begin{equation}
r([\phi]) = [\phi|_H] \label{restrict}
\end{equation}
Here, $H^1(H,\BZ G)$ is the quotient of the set of all functions $\phi \co H \to \BZ G$ satisfying $\phi(hh')= \phi(h) + h\phi(h')$ by the subset of such $\phi$ that are of the form $\phi(h)=(h-1)y$ for some $y\in \BZ G$.

\begin{theorem}[Swarup \cite{swarup1977relative}]\label{swarup}
Let $G$ be a finitely generated group, and let $H_1,\ldots,H_m$ be subgroups of $G$. 
For each $j$, let
$$r_j \co H^1(G,\BZ G) \to H^1(H_j, \BZ G)$$
denote the restriction map. If the intersection $\bigcap_j \mathrm{ker}(r_j) \neq 1$, then there is an action of $G$ on a simplicial tree $T$, with finite edge stabilizers, no edge inversions, no global fixed point, and such that each $H_j$ is contained in the stabilizer of some vertex $p_j \in T$.
\end{theorem}

Equivalently, as stated in \cite{swarup1977relative}, the conclusion is that $G \cong A \star_C B$ or $G\cong A \star_C$ in such a way that each $H_j$ is conjugate into $A$ or $B$. 
The following proposition interprets Swarup's theorem in terms of end spaces, at least when the subgroups are finitely generated.

\begin{proposition}\label{swarup in english}
Let $H\subset G$ be finitely generated, and let $r$ be the restriction map of \eqref{restrict}. Under the isomorphism $H^1(G,\BZ G) \to C(\CE(G),\BZ)/Z \) given by \( [\phi_x] \longmapsto [\bar x]$ defined in Theorem \ref{bigthm}, the kernel of $r$ maps to the subset of $C(\CE(G),\BZ)/Z$ consisting of all functions $f : \CE(G) \to \BZ$ such that for each right coset $Hg$, the function $f$ is constant on  $\overline{H g} \cap \CE(G)$.
\end{proposition}

Here, $\overline{Hg}$ is the closure of \( Hg \)  in $G \cup \CE(G)$. 
Given multiple subgroups $H_j$ as in Swarup's theorem, the intersection of the kernels is the set of all $f$ that are constant on $\overline{H_jg} \cap \CE(G)$ for each $j$ and each right coset $H_j g$.

\begin{proof}[Proof of Proposition~\ref{swarup in english}]
By Lemma~\ref{cocycles}, any $1$-cocycle with \( \BZ G \)-coefficients is of the form $\phi_x$, where $x\in \BZ^G$ has finite boundary. 
The map $r$ takes $[\phi_x]$ to the class in $H^1(H,\BZ G)$ of the restriction
$\phi_x|_H : H \to \BZ G$ given by  
\[ \phi_x|_H(h)=(h-1), \]
which is trivial in $H^1(H_j ,\BZ G)$ exactly when there is an element $y\in \BZ G$ such that $\phi_x(h)=\phi_y(h)$ for all $h\in H$. This last equality is equivalent to saying that $h(x-y) = x-y$ for all $h\in H$, or saying that $x-y $ is invariant under the left action of $H$.
So, we have that $r ([\phi_x])=0$ if and only if  the function $x$ is left $H$-invariant outside of a finite subset of $G$. 

If $x$ is left $H$-invariant outside a finite subset of $G$, then on each right coset $Hg$, the function $x$ is constant outside of a finite subset, and hence the boundary map $\bar x$ is constant on $\overline {Hg} \cap \CE(G)$. 

Conversely, suppose $\bar x$ is constant on $\overline {Hg} \cap \CE(G)$ for each $g \in G$. 
Fix a finite generating set \( S \) for \( G \). 
As $H$ is finitely generated, we can pick $m \in \BN$ large enough so that \( H \) is generated by elements of $H\subset G$ with word length at most $m$ in the generating set \( S \) of \( G \). Note that any two points in a coset $Hg$ can then be joined in $Cay(G)$ by a path that stays in the $m$-neighborhood of $Hg$.
Call a coset $Hg$ \emph{relevant} if it intersects the $m$-neighborhood $\CN_m(\partial x)$ of $\partial x$, and call it \emph{irrevelant} otherwise. 
As \( \partial x \) is finite, there are only finitely many relevant cosets; moreover, each irrelevant coset lies in a single component of $Cay(G) \smallsetminus \CN_m(\partial x)$.
Let $F \subset G$ be the union of $\CN_m(\partial x)$  and all intersections $Hg \cap U$ that are finite sets, where $Hg$ is any relevant coset and $U$ is any component of $Cay(G) \smallsetminus \CN_m(\partial x)$. 
Since there are only finitely many $U$ and finitely many relevant $Hg$, the set $F$ is finite. 

We claim that $x$ is left $H$-invariant on the complement of $F$, meaning that for every right coset $Hg$, the function $x$ is constant on $Hg \smallsetminus F$. 
If $Hg$ is irrelevant, then it lies in a single component of $Cay(G) \smallsetminus \CN_m(\partial x)$, and hence any two points in $Hg$ can be joined by a path avoiding $\partial x$, implying that $x$ is constant on $Hg$ as desired. 
If $Hg $ is relevant, then since $\bar x$ is constant on $\overline {Hg} \cap \CE(G)$,  it follows that on all components of $Cay(G) \smallsetminus \CN_m(\partial x)$ that $Hg$ intersects in an infinite set, the value of $x$ is the same constant. Hence $x$ is constant on $Hg \smallsetminus F$ as desired.
\end{proof}

Combining Swarup's theorem with Proposition~\ref{swarup in english}, we get Theorem \ref{Thm:relativestallings}.

  \bibliographystyle{amsplain}
  \bibliography{total}
  \end{document}